\documentclass[a4paper,oneside]{amsart}
\usepackage{geometry}

\usepackage{mathtext}
\usepackage[T1]{fontenc}
\usepackage{color}
\usepackage{amsmath}
\usepackage{amsfonts,amssymb,amsthm}
\usepackage[pdftex]{graphicx}
\usepackage{comment}
\usepackage{float}
\usepackage{todonotes}
\usepackage{microtype}
\usepackage{mathdots}
\usepackage{verbatim}
\usepackage{ytableau}

\usepackage[colorlinks=true, linkcolor=blue, citecolor=blue, linktoc=page]{hyperref}

\usepackage{tikz}
\usetikzlibrary{decorations.markings}
\tikzset{middlearrow/.style={
		decoration={markings,
			mark= at position 0.5 with {\arrow{#1}} ,
		},
		postaction={decorate}
	}
}
\tikzset{firstthirdarrow/.style={
		decoration={markings,
			mark= at position 0.33 with {\arrow{#1}} ,
		},
		postaction={decorate}
	}
}
\tikzset{secondthirdarrow/.style={
		decoration={markings,
			mark= at position 0.66 with {\arrow{#1}} ,
		},
		postaction={decorate}
	}
}
\usetikzlibrary{cd}

\usepackage{enumerate}

\DeclareMathOperator{\Span}{Span}

\DeclareMathOperator{\Sp}{Sp}

\DeclareMathOperator{\SL}{SL}

\DeclareMathOperator{\SO}{SO}

\DeclareMathOperator{\Ein}{Ein}
\DeclareMathOperator{\Mink}{Mink}
\DeclareMathOperator{\Iso}{Iso}
\renewcommand{\det}{\mathrm{det}}

\usepackage{color}
\definecolor{NoteColor}{rgb}{1,0,0}

\newcommand{\R}{\mathbb R}

\newcommand{\RR}{\mathbb R}
\renewcommand{\SS}{\mathbb S}

\theoremstyle{plain}
\newtheorem{teo}{Theorem}[section]
\newtheorem{theorem}[teo]{Theorem}

\newtheorem{corollary}[teo]{Corollary}

\newtheorem{lemma}[teo]{Lemma}
\newtheorem{prop}[teo]{Proposition}
\newtheorem{proposition}[teo]{Proposition}

\newtheorem{question}[teo]{Question}

\newtheorem{teoletter}{Theorem}
\newtheorem*{questionletter}{Question}

\theoremstyle{definition}

\newtheorem{definition}[teo]{Definition}
\newtheorem{ex}{Example}[section]

\theoremstyle{remark}
\newtheorem{rem}[teo]{Remark}
\newtheorem{remark}[teo]{Remark}

\usepackage{xparse}

\renewcommand{\bar}{\overline}

\usepackage[
backend=biber,
giveninits=true,
style=alphabetic,
doi=false,isbn=false,url=false,eprint=false
]{biblatex}
\begin{document}
	
	\title{Connected components of the space of flags of $\SO_0(p,q)$ transverse to a fixed pair and restrictions on Anosov subgroups}

	\author[C. Kineider]{Clarence Kineider}
	
	\author[R. Troubat]{Roméo Troubat}
	
	
	
	\address{Clarence Kineider \\ Max Planck Institute for Mathematics in the Sciences, Leipzig, Germany}
	\email{kineider@mis.mpg.de}
	
	\address{Roméo Troubat \\ Institut des Hautes \'Etudes Scientifiques, Bures-Sur-Yvette, France}
	\email{troubat@ihes.fr}
	
	\begin{abstract}
        We count and give a parametrization of connected components in the space of flags transverse to a given transverse pair in every flag varieties of $\SO_0(p,q)$. We compute the effect the involution of the unipotent radical has on those components and, using methods of Dey--Greenberg--Riestenberg, we show that for certain parabolic subgroups $P_{\Theta}$, any $P_{\Theta}$-Anosov subgroup is virtually isomorphic to either a surface group of a free group. We give examples of Anosov subgroups which are neither free nor surface groups for some sets of roots which do not fall under the previous results. As a consequence of the methods developed here, we get an explicit computation of some Pl\"ucker coordinates to check if a unipotent matrix in $\SO_0(p,q)$ belong to the $\Theta$-positive semigroup $U_\Theta^{>0}$ when $p\neq q$.
	\end{abstract}

	\thanks{This project has received funding from the European Research Council (ERC) under the European Union’s Horizon 2020 research and innovation programme (grant agreement No 101018839), as well as from the \'Ecole Normale Supérieure and the Strasbourg Université.}

	\maketitle
	
	\tableofcontents
 
	\section*{Introduction}

    In his celebrated paper \cite{Lab04}, Labourie introduced Anosov representations of a surface group as a way to generalize Hitchin representations, introduced by Hitchin in \cite{Hit92}. The notion of Anosov representations was then expanded to any hyperbolic group by Guichard--Wienhard in \cite{GW12}. For a finitely generated hyperbolic group $\Gamma$ and a semi-simple Lie group $G$, when a representation $\rho : \Gamma\to G$ is Anosov with respect to a parabolic subgroup $P\subset G$ (or $P$-Anosov), there exists a $\rho$-equivariant \emph{limit map} $\xi: \partial_\infty\Gamma\to G/P$ from the Gromov boundary of $\Gamma$ to the flag variety $G/P$. This limit map has  a \emph{transversality} property: given two distinct points in $\partial_\infty\Gamma$, their image are transverse in $G/P$ (when $P$ is self-opposite). This imply that Anosov-representations are discrete and faithful. When the parabolic subgroup is the Borel subgroup (i.e. the parabolic subgroup of $G$ associated to the whole set of simple roots), the representation is usually called Borel-Anosov and its image is called a Borel-Anosov subgroup. This is the most restrictive Anosov property, as a Borel-Anosov representation is in particular $P$-Anosov for any parabolic subgroup $P$. Known examples of Borel-Anosov representations include representations of surface groups and of free groups, for instance through the inclusion of co-compact lattices and Schottky groups of $\SL(2, \R)$ via the irreducible representation $\SL(2, \R) \rightarrow \SL(d,\R)$. A question was then asked by Andrés Sambarino regarding the abstract structure of groups admitting Borel-Anosov representations in $\SL(d, \R)$.

    \begin{questionletter}
        Let $\Gamma$ be a Borel-Anosov subgroup of $\SL(d,\R)$. Is $\Gamma$ virtually isomorphic to either a free group or a surface group?
    \end{questionletter}

    This question was answered positively for certain values of $d$, first by Canary--Tsouvalas and Tsouvalas in \cite{Canary_2020}, \cite{tsouvalas2019borelanosovrepresentationsdimensions} for $d=3, 4$ and $d=2$ mod $4$, then by Dey in \cite{dey2024borelanosovsubgroupsrm} for $d \neq 5$ and $d \neq \pm 1$ mod $8$.

    For $F_0$ and $F_{\infty}$ two transverse elements in the space of full flags $\SL(d, \R)/B$, where $B$ is the Borel subgroup of $\SL(d, \R)$, the space $\Omega(F_{\infty})$ of flags that are transverse to $F_\infty$ is parametrized by the unipotent radical $U$ of $B$ via the map $g \in U \mapsto g \cdot F_0$. This parametrization induces a natural involution defined by $i(g \cdot F_0) = g^{-1} \cdot F_0$ which preserves transversality with $F_0$, thus inducing an involution $i : \Omega(F_0) \cap \Omega(F_{\infty})\to   \Omega(F_0) \cap \Omega(F_{\infty})$. In \cite{dey2024borelanosovsubgroupsrm}, Dey uses this involution to obtain an upper bound on the dimension of boundaries of Borel-Anosov subgroups of $\SL(d, \R)$. In \cite{dey2023restrictionsanosovsubgroupssp2nr}, Dey--Greenberg--Riestenberg generalize this method to any pair $(G, P_{\Theta})$ where $G$ is a semi-simple Lie group and $P_{\Theta} \subset G$ is a parabolic subgroup such that $G/P_{\Theta}$ is self-opposite, thus obtaining the following theorem:

    \begin{theorem}[Dey--Greenberg--Riestenberg, \cite{dey2023restrictionsanosovsubgroupssp2nr}]\label{DGR}
        Assume that $$i : \pi_0(\Omega(F_0) \cap \Omega(F_{\infty})) \rightarrow \pi_0(\Omega(F_0) \cap \Omega(F_{\infty}))$$ does not have any fixed point. Then any $P_{\Theta}$-Anosov subgroup of $G$ is virtually isomorphic to either a free group or a surface group.
    \end{theorem}

    They then apply this result to $G = \Sp(2n,\mathbb{R})$ obtaining that any $P_{\Theta}$-Anosov subgroup of $\Sp(2n, \R)$ is virtually isomorphic to either a surface group or a free group when $\Theta$ contains an odd root. The idea behind this theorem is that given two points $F_0$ and $F_\infty$ in the image of the limit map of a $P_\Theta$-Anosov representation, the complement of these two points in the image is either contained in a single connected component of $\Omega(F_0) \cap \Omega(F_{\infty})$ stable under the involution $i$ or in a pair of connected components of $\Omega(F_0) \cap \Omega(F_{\infty})$ swapped by $i$. In the latter case, topological arguments show that the image of the limit map have to be contained inside a topological circle. A classical result about Gromov-hyperbolic groups imply that a finitely generated Gromov hyperbolic group whose boundary is contained in a circle is virtually isomorphic either to a surface group or a free group.
    
    There has recently been a surge of interest in the groups $\SO(p,q)$, first due to the discovery of $\Theta$-positive representations of a surface group into $\SO(p,q)$ by Guichard--Wienhard in \cite{GW24}, and then by the discovery, first by Barbot for $q=1$ in \cite{Barbot} and by Beyrer-Kassel for $q > 1$ in \cite{BK25}, of so-called "higher higher Teichmuller spaces" of $\SO(p,q)$, meaning connected components of the caracter variety of a hyperbolic group of cohomological dimension greater than $2$ into $\SO(p,q)$ made entirely of discrete and faithful representations. Those higher higher Teichmuller spaces are exactly the $P_1$-Anosov representations of a group with no infinite nilpotent normal subgroup of virtual cohomological dimension $p$ into $\SO(p,q)$. 
    
    Our goal in this paper will be to apply Theorem \ref{DGR} whenever possible when $G = \SO_0(p,q)$ (for $p\geq q\geq 2$) is the connected component of the identity in the group $\SO(p,q)$ of isometries preserving a quadratic form of signature $(p,q)$, and when $P_{\Theta}$ is any parabolic subgroup of $\SO_0(p,q)$ associated to a subset $\Theta$ of the simple restricted roots. The first step in this endeavor will be to count, parameterise and overall get a better understanding of the connected components of $\Omega(F_0) \cap \Omega(F_{\infty})$, the space of flags transverse to a fixed pair of transverse flags $(F_0,F_\infty)$. The second step will be to compute the action of the involution $i$ on those connected components. While the computation of the involution does not always allow oneself to conclude that any corresponding Anosov subgroups are virtually isomorphic to surface groups or free groups, one may still get a similar result if they first manage to show that the boundary of the group is contained within a pair of connected components which are exchanged by the involution. In what follows, we refer to Sections \ref{sec:flag_var_pq} and \ref{sec:flag_var_qq} for the conventions regarding the labeling of the roots and parabolic subgroups of $\SO_0(p,q)$.

    \begin{teoletter}\label{th:CC_pq}
        Assume that~$p > q+1$ and $q \geq 2$. Then the number of connected components of~$\Omega(F_0) \cap \Omega(F_{\infty})$ and the action of~$i$ are as follows:
       \begin{itemize}
         \item If~$\Theta$ does not contain the last root and is not equal to~$\{1,\dots,q-1\}$, there are~$2^{|\Theta|}$ connected components, all of which are stable by~$i$.
         \item If~$\Theta$ contains the last root and is not equal to~$\{1,\dots,q\}$, there are~$2^{|\Theta|-1}$ connected components, all of which are stable by~$i$.
         \item If~$q > 3$,~$\Theta = \{1,\dots,q-1\}$ or~$\Theta = \{1,\dots,q\}$, there are~$3 \times 2^{q-1}$ connected components. Of them, $2^{q-1}$ are stable by~$i$ while the rest are not stable.
         \item If~$q=3$,~$\Theta = \{1,2\}$ or~$\Theta = \{1,2,3\}$, there are $11$ connected components. Six of them are stable by~$i$ while the rest are not stable. 
         \item If~$q=2$,~$\Theta=\{1\}$ or~$\Theta = \{1,2\}$, there are three connected components. One is stable by~$i$ and the other two are not. 
       \end{itemize}
    \end{teoletter}
    
    It follows that when~$p>q+1$, there is always at least one connected component of~$\Omega(F_0) \cap \Omega(F_{\infty})$ which is stable by~$i$, thus the methods of Dey--Greenberg--Riestenberg do not yield general results on the structure of any Anosov subgroups.
    
    \begin{teoletter}\label{th:CC_qq+1}
        Assume that~$p = q+1$ and $q \geq 2$. Then the number of connected components of~$\Omega(F_0) \cap \Omega(F_{\infty})$ and the action of~$i$ are as follows:
       \begin{itemize}
         \item If~$\Theta$ is not equal to~$\{1,\dots,q-1\}$ nor to~$\{1,\dots,q\}$, there are~$2^{|\Theta|}$ connected components. If $\Theta$ does not contain $q$, all of them are stable by $i$. In the case where $\Theta$ contains $q$, if $q = 0$ or $3$ mod $4$, all connected components are stable by $i$ while if $q = 1$ or $2$ mod $4$, none of them are.
         \item If~$q > 3$,~$\Theta = \{1,\dots,q-1\}$, there are~$3 \times 2^{q-1}$ connected components. Half of the components are stable by~$i$ while the other half are not stable.
         \item If~$q > 3$,~$\Theta = \{1,\dots,q\}$, there are~$(q+5)2^{q-1}$ connected components. When~$q=1$ or~$2$ mod~$4$, none of the components are stable by~$i$. When~$q=3$ or~$0$ mod~$4$,~$2^q$ of the components are stable by~$i$ and~$(q+3)2^{q-1}$ are not.
         \item If~$q=3$,~$\Theta = \{1,2\}$, there are $11$ connected components. Six of the components are stable by~$i$ while the rest are not stable.
         \item If~$q=3$,~$\Theta=\{1,2,3\}$, there are~$30$ connected components. Six of the components are stable while the rest are not stable.
         \item If~$q=2$,~$\Theta=\{1\}$, there are three connected components. One of the components is stable by~$i$ and the other two are not.
         \item If~$q=2$,~$\Theta = \{1,2\}$, there are eight connected components. None of the components are stable by~$i$.
       \end{itemize}
    \end{teoletter}

    The group $\SO(q+1,q)$ is split and when $\Theta = \{1,\dots,q\}$, $P_\Theta$ is its Borel subgroup. In that case, the number of connected components in $\Omega(F_0) \cap \Omega(F_{\infty})$ was already known from work of Gekhtman--Shapiro--Vainshtein in \cite{GSV03}. 
    
    \begin{teoletter}\label{th:CC_qq}
        Assume that~$p =q$ and $q \geq 2$. Then the number of connected components of~$\Omega(F_0) \cap \Omega(F_{\infty})$ is as follows:
       \begin{itemize}
         \item If $\Theta$ does not contain any of the last two roots, there are $2^{|\Theta|}$ connected components, all of which are stable by $i$.
         \item If $\Theta$ contains only one of the last two roots, $\SO(q,q)/P_{\Theta}$ is self-opposite if and only if $q$ is even. In this case, there are $2^{|\Theta|}$ connected components. When $q = 0$ mod $4$, all of them are stable by $i$ whereas when $q = 2$ mod $4$, none of them are.
         \item If $\Theta$ contains both of the last two roots and is not equal to $\{1,\dots,q-2,q,q'\}$, there are $2^{|\Theta|}$ connected components. 
         \item If $q > 3$ and $\Theta$ is equal to $\{1, \dots, q-2, q, q'\}$, there are $3 \times 2^q$ connected components.
         \item If $q=3$ and $\Theta$ is equal to $\{1,3,3'\}$, there are $20$ connected components.
         \item If $q=2$ and $\Theta$ is equal to $\{2, 2'\}$, there are $4$ connected components. None of the components are stable by $i$.
       \end{itemize}
    \end{teoletter}
    
    When it is not specified, we did not compute the involution on the present work. However in those cases, recent results by Evans--Riestenberg \cite{evans2025transversespheresflagmanifolds} show that at least one connected component must be stable. The group $\SO(q,q)$ is split, and the component count in the case when $\Theta = \{1,\dots,q-2,q,q'\}$ and $P_{\Theta}$ is the Borel subgroup was achieved by Zelevinsky \cite{Zel00}. Though we do not do it explicitly here, our methods would yield the same result.

    Combining those results with theorem \ref{DGR}, we obtain the following theorems:

    \begin{teoletter}\label{rigid1}
        Let $q \geq 2$. Assume that $q = 1\text{ mod }4$ or $q=2\text{ mod }4$ and let $\Gamma$ be a $P_{\Theta}$-Anosov subgroup of $\SO_0(q+1,q)$ with $\Theta$ containing the last root. Then $\Gamma$ is virtually isomorphic to either a free group or a surface group.
    \end{teoletter}

    \begin{teoletter}\label{rigid2}
        Let $q \geq 2$. Assume that $q = 2$ mod $4$ and let $\Gamma$ be a $P_{\Theta}$-Anosov subgroup of $\SO_0(q,q)$ with $\Theta$ containing one of the last two roots. Then $\Gamma$ is virtually isomorphic to either a free group or a surface group.
    \end{teoletter}

    When taking $P_{\Theta}$ to be the Borel subgroup of $\SO_0(q+1,q)$, the theorem \ref{rigid1} is implied by the results in \cite{dey2024borelanosovsubgroupsrm} via the inclusions $\SO_0(q+1,q) \subset \SL(2q+1, \R)$. However when taking $P_{\Theta}$ to be the Borel subgroup of $\SO_0(q,q)$, the theorem \ref{rigid2} give a stronger result than the one obtained by applying the result of \cite{dey2024borelanosovsubgroupsrm} to the embedding $\SO_0(q,q) \subset \SL(2q, \R)$ as a Borel-Anosov subgroup of $\SO_0(q,q)$ may not be Borel-Anosov in $\SL(2q, \R)$.

    We show that to determine the connected component of $\Omega(F_0) \cap \Omega(F_{\infty})$ to which a given flag $g\cdot F_0$ belong to, one only has to check the signs of some explicit minors (which are Pl\"ucker coordinates) of the matrix $g$. In particular, we give an explicit set of signs these minors have to be in order for $g$ to be a $\Theta$-positive matrix in the sense of Guichard--Wienhard in \cite{GW24} (see Remark \ref{rem:Theta_positive_cells}). This yields a parametrization of the $\Theta$-positive semigroup $U_\Theta^{>0}$ that is different from the one given in \cite{GW24}.

    Finally, we give examples of $P_{\Theta}$-Anosov subgroups in $\SO_0(p,q)$ for $\Theta = \{1,...,q-2\}$ and $\Theta = \{1,...,q\}$ which are not virtually isomorphic either to a surface group nor to a free group.

    \begin{teoletter}\label{exanosov}
        Let $p > q \geq 3$. There exists a $P_{1, ...,q-2}$-Anosov subgroup of $\SO_0(p,q)$ which is isomorphic to the free product of a surface group by an infinite cyclic group.
    \end{teoletter}

    \begin{teoletter}\label{exanosov2}
        When $q \geq 2$ and $p \geq 3q+1$, there exists a $P_{1, ...,q}$-Anosov subgroup of $\SO_0(p,q)$ which is isomorphic to the free product of a surface group by an infinite cyclic group.
    \end{teoletter}

    In Section 1, we define and give a parametrization of the various flag varieties associated to $\SO_0(p,q)$, as well as recall basic terminology of the geometry of the pseudo-Riemannian Einstein universe. In Section 2, we establish the equations characterizing transversality to a given point in an affine chart of any flag variety, both via a geometric approach and a computational approach. In Section 3, we count the number of connected components in $\Omega(F_0) \cap \Omega(F_{\infty})$ for every flag variety associated to $\SO_0(p,q)$, thus partially proving Theorems \ref{th:CC_pq}, \ref{th:CC_qq+1} and \ref{th:CC_qq}. In Section 4, we compute the action of the involution on those connected components, thus finishing the proofs of the previous theorems and establishing Theorems \ref{rigid1} and \ref{rigid2}. In Section 5, we construct the Anosov subgroups described in Theorems \ref{exanosov} and \ref{exanosov2}.

    \section*{Acknowledgments}

    We would like to thank Olivier Guichard, Subhadip Dey and Max Riestenberg for their helpful comments during the preparation of this article. We are also grateful towards Olivier Schiffmann for organizing the math camp during which most of the ideas presented here arose.

    \section{Flag varieties of \texorpdfstring{$\SO_0(p,q)$}{SO(p,q)}} 

    \subsection{Conventions and notations}

    Let $p\geq q\geq 1$ and let $\mathcal{B} = (e_1,\dots,e_{q},x_{q+1},\dots,x_{p-q},\Tilde{e}_{q},\dots,\Tilde{e}_1)$ be a basis of $\RR^{p+q}$.
    Let $Q_{p,q}$ be the non-degenerate quadratic form on $\RR^{p+q}$ of signature $(p,q)$ which matrix in the basis $\mathcal{B}$ is $$M_{p,q} = \begin{pmatrix}
        0 & 0 & J\\
        0 & I_{p-q} & 0\\
        J & 0 & 0
    \end{pmatrix}$$ 
    where $J$ is the following $q$-by-$q$ square matrix: $$ J = \begin{pmatrix}
        0 & \cdots & 0 & 1\\
        \vdots & \iddots & 1 & 0\\
        0 & \iddots & \iddots & \vdots\\ 
        1 & 0 & \dots & 0
    \end{pmatrix}.$$
    When the context is clear we will write $Q$ instead of $Q_{p,q}$, we will denote by $B_{p,q}$ (or just $B$) the associated bilinear form, and for any subspace $F$ of $\RR^{p+q}$, we will denote by $F^\perp$ its orthogonal:
    $$F^\perp = \left\lbrace v\in \RR^{p,q} ~\vert~ \forall w\in F, B(v,w) = 0\right\rbrace.$$
    A subspace $F$ will be called \emph{isotropic} if $Q|_F = 0$. In particular, note that $\Span(e_1,\dots e_q)$ is a maximal isotropic subspace and that the quadratic form $Q$ restricted to $\Span(x_{q+1},\dots,x_{p-q})$ is positive definite. 
    
    We denote by $\RR^{p,q}$ the space $\RR^{p+q}$ endowed with the quadratic form $Q_{p,q}$. The group $\SO(p,q)$ is the group of orientation-preserving isometries of $\RR^{p,q}$, and $\SO_0(p,q)$ will denote the identity component of $\SO(p,q)$. The Lie algebra of $\SO_0(p,q)$ is the space 
    $$\mathfrak{so}(p,q) = \left\lbrace A\in\mathcal{M}_{p+q}(\RR) ~\vert~ A^TM_{p,q}+M_{p,q}A = 0\right\rbrace.$$
    For any non-negative integer $k\geq 0$, an isotropic subspace of $\RR^{p,q}$ of dimension $k+1$ will be called a $k$-\emph{photon}. Sometime a $1$-photon will simply be called a \emph{photon}, and a $0$-photon will be called a \emph{point}. A vector $v\in\RR^{p,q}$ will be called \emph{timelike} if $Q(v)<0$, \emph{spacelike} if $Q(v)>0$ and \emph{lightlike} if $Q(v)=0$. The set of all lightlike vectors form a cone, called the \emph{lightcone} (of the origin), denoted by $C(0)$. Given $v\in\RR^{p,q}$, the lightcone of $v$ is the translation of $C(0)$:
    $$C(v) = \left\lbrace w\in \RR^{p,q} ~\vert~ Q(w-v) = 0\right\rbrace.$$
    The lightcone $C(0)$ divides $\RR^{p,q}$ into two parts: the set of timelike vectors and the set of spacelike vectors. The set of spacelike vectors is a connected component of $\RR^{p,q}\backslash C(0)$ when $p > 1$ and is the union of two connected components when $p=1$. Likewise, the set of timelike vectors is a connected component of $\RR^{p,q}\backslash C(0)$ if $q>1$ and has two connected components if $q=1$. When the space of timelike vectors has two connected components, we will call a vector in the component containing $e_1-\Tilde{e}_1$ of \emph{future type} and a vector in the component of $\Tilde{e}_1-e_1$ of \emph{past type}.
    
    Depending on the relative values of $p$ and $q$, the group $\SO_0(p,q)$ may have different properties. We will assume in the following that $q \geq 2$, the case $q=1$ being already well studied. Since we want to cover every non-compact group of this family, we will often have to do a case by case study, depending on whether $p=q$, $p=q+1$ or $p>q+1$. When $p=q$ or $p=q+1$, the group $\SO_0(p,q)$ is split-real. When $p=q+1$, $\SO_0(p,q)$ has a (restricted) root system of type $B_q$ (see Figure \ref{fig:Bn_diagram}), and when $p=q$ the group $\SO_0(q,q)$ has a root system of type $D_q$ (see Figure \ref{fig:Dn_diagram}). Because of the different (restricted) root systems, the flag varieties associated to $\SO_0(p,q)$ will be slightly different whether $p=q$ or $p\neq q$. 
    
    \subsubsection{Flag varieties when \texorpdfstring{$p\neq q$}{p =! q}}\label{sec:flag_var_pq}
    
\begin{figure}[h!]
    \centering 
      \includegraphics[width=.35\linewidth]{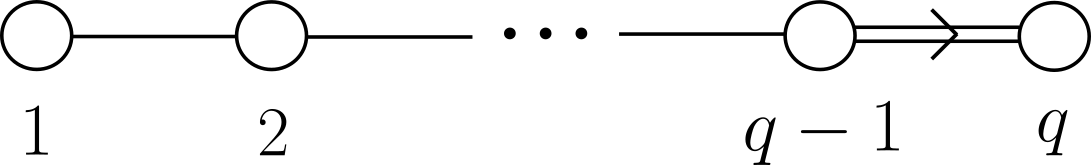}
      \caption{The Dynkin diagram $B_q$ associated to $\SO(p,q)$ for $p\geq q+1$.}\label{fig:Bn_diagram}
    \end{figure}

    We start by describing the different flag varieties of $\SO_0(p,q)$ when $p\geq q+1$. The Dynkin diagram corresponding to the (restricted) root system is $B_q$, with the set of simple roots $\Delta$ labeled by integers from 1 to $q$ as in Figure \ref{fig:Bn_diagram}.
    
    \begin{definition}
        Let $q \geq 2$ and $p\geq q+1$. Let $k\geq 1$ be a positive integer and let $\Theta  = (i_1,\dots,i_k)$ be a $k$-tuple of integers satisfying $1\leq i_1 <\dots <i_k\leq q$. A $\Theta $-flag in $\RR^{p,q}$ is a $k$-tuple $F = (F^{i_1},\dots,F^{i_k})$ such that for all $1\leq \ell\leq k$, $F^{i_\ell}$ is a $(i_\ell-1)$-photon (i.e a $i_{\ell}$-dimensional isotropic subspace of $\R^{p,q}$) and for all $1\leq \ell\leq k-1$, $F^{i_\ell}\subset F^{i_{\ell+1}}$. The space of all $\Theta $-flags of $\RR^{p,q}$, called the $\Theta $-flag variety of $\SO_0(p,q)$, will be denoted by $\mathcal{F}^\Theta $. When $\Theta  = (1,2,\dots,q)$, we will call $\mathcal{F}^\Theta $ the \emph{full flag variety}.
    \end{definition}

    Let $\Theta  = (i_1,\dots,i_k)$, we define the two \emph{standard} $\Theta $-flags $F^\Theta _0$ and $F^\Theta _\infty$, or just $F_0$ and $F_\infty$ when the context is clear, as follows:
    $$ F^\Theta _0 = (F_0^{i_\ell})_{1\leq \ell \leq k},\text{ where } F_0^{i_\ell} = \Span(\Tilde{e}_1,\dots,\Tilde{e}_{i_\ell})$$
    and
    $$ F^\Theta _\infty = (F_\infty^{i_\ell})_{1\leq \ell \leq k},\text{ where } F_\infty^{i_\ell} = \Span(e_1,\dots,e_{i_\ell}).$$

    The group $\SO_0(p,q)$ acts smoothly and transitively on $\mathcal{F}^\Theta $, and the stabilizer of $F^\Theta _0$ (resp. $F_\infty^\Theta $) is a parabolic subgroup denoted by $P_\Theta $ (resp. $P_\Theta ^{opp}$). This parabolic subgroup is the one associated to the subset of simple roots $\Theta \subset\Delta$. Thus, the space $\mathcal{F}^\Theta $ is diffeomorphic to $\SO_0(p,q)/P_\Theta $ which is a smooth compact manifold. When $\Theta =(1,2,\dots,q)$, $P_\Theta $ is a minimal parabolic subgroup of $\SO_0(p,q)$, which is called a \emph{Borel subgroup} when $p=q+1$. When $p\neq q$, every parabolic subgroup of $\SO_0(p,q)$ is self-opposite (i.e. $P_\Theta $ is conjugated to $P_\Theta ^{opp}$).
    
    \subsubsection{Flag varieties when \texorpdfstring{$p=q$}{p=q}}\label{sec:flag_var_qq}
    
    When $p=q\geq 2$, the group $\SO_0(q,q)$ does not act transitively on the set of $(q-1)$-photons, but has two distinct orbits. We refer to Section \ref{sec:einstein} for a geometric description of these two orbits. The stabilizers of a $(q-1)$-photon in either of these two orbits are both parabolic subgroups of $SO_0(q,q)$ associated to the two forking roots in the Dynkin diagram $D_q$, which we will hence label $q$ and $q'$, the usual convention being to label one of these by $q-1$, but we want to keep a correspondence between the label of the root defining the parabolic subgroup and the dimension of the isotropic subspace stabilized. The labeling of the simple roots $\Delta$ of $\SO_0(q,q)$ is then:
    
    \begin{figure}[h!]
    \centering 
      \includegraphics[width=.35\linewidth]{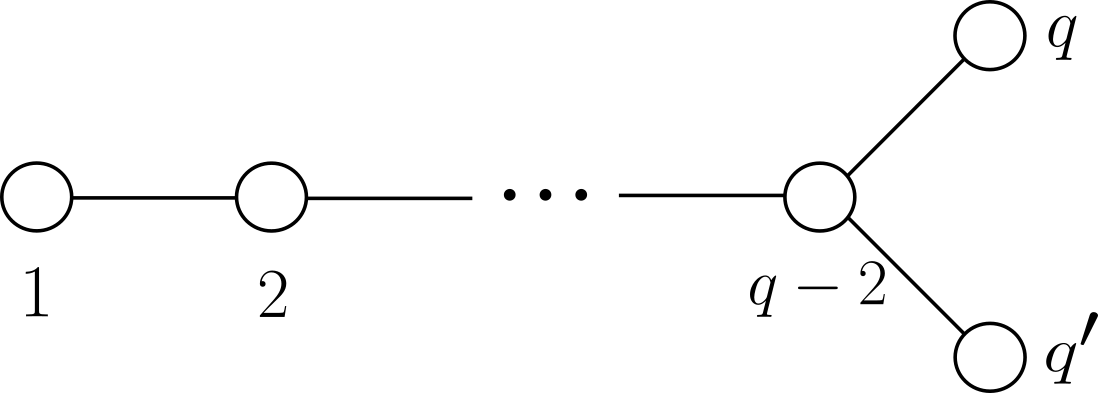}
      \caption{The Dynkin diagram $D_q$ associated to $\SO(q,q)$ with our convention of root labeling.}\label{fig:Dn_diagram}
    \end{figure}
    
    \begin{definition}
        Let $q\geq 2$. Let $k\geq 1$ be a positive integer and let $\Theta  = (i_1,\dots,i_k)$ be a $k$-tuple of distinct elements in $\lbrace 1,\dots, q-2,q,q'\rbrace$ satisfying $1\leq i_1 <\dots <i_k\leq q$ (with the convention that $q-2<q'<q$). A $\Theta $-flag in $\RR^{q,q}$ is a $k$-tuple $F = (F^{i_1},\dots,F^{i_k})$ such that for all $1\leq \ell\leq k$, $F^{i_\ell}$ is a $(i_\ell-1)$-photon (or a $(q-1)$-photon if $i_\ell=q'$) and for all $1\leq \ell\leq k-1$, $F^{i_\ell}\subset F^{i_{\ell+1}}$ (unless $i_\ell=q'$ and $i_{\ell+1}=q$). The space of all $\Theta $-flags of $\RR^{q,q}$, called the $\Theta $-flag variety of $\SO_0(q,q)$, will be denoted by $\mathcal{F}^\Theta $. When $\Theta  = (1,2,\dots,q-2,q',q)$, we will call $\mathcal{F}^\Theta $ the \emph{full flag variety}.
    \end{definition}
    
    Let $\Theta  = (i_1,\dots,i_k)$, we define the two \emph{standard} $\Theta $-flags $F^\Theta _0$ and $F^\Theta _\infty$, or just $F_0$ and $F_\infty$ when the context is clear, as follows:
    $$ F^\Theta _0 = (F_0^{i_\ell})_{1\leq \ell \leq k},\text{ where } F_0^{i_\ell} = \left\lbrace\begin{array}{l}
         \Span(\Tilde{e}_1,\dots,\Tilde{e}_{i_\ell}) \text{ if } i_\ell\neq q'\\
         \Span(\Tilde{e}_1,\dots,\Tilde{e}_{q-1}, e_q) \text{ if } i_\ell = q'
    \end{array}\right.$$
    and
    $$ F^\Theta _\infty = (F_\infty^{i_\ell})_{1\leq \ell \leq k},\text{ where } F_\infty^{i_\ell} = \left\lbrace\begin{array}{l}
         \Span(e_1,\dots,e_{i_\ell}) \text{ if } i_\ell\neq q'\\
         \Span(e_1,\dots,e_{q-1}, \Tilde{e}_q) \text{ if } i_\ell = q'
    \end{array}\right. .$$

    The group $\SO_0(q,q)$ acts smoothly and transitively on $\mathcal{F}^\Theta $, and the stabilizer of $F^\Theta _0$ (resp. $F_\infty^\Theta $) is a parabolic subgroup denoted by $P_\Theta $ (resp. $P_\Theta ^{opp}$). Thus, the space $\mathcal{F}^\Theta $ is diffeomorphic to $\SO_0(q,q)/P_\Theta $ which is a smooth compact manifold. When $\Theta =(1,2,\dots,q-2,q',q)$, $P_\Theta $ is a Borel subgroup of $\SO_0(q,q)$. 
    
    A consequence of this is that the Grassmannian of isotropic $(q-1)$-dimensional subspaces (i.e. $(q-2)$-photons) does not arise immediately as a flag variety of $\SO_0(q,q)$. In fact it is one of the flag varieties of $\SO_0(q,q)$, but not one associated to a single root:
    
    \begin{prop}\label{prop:sub-max_photon}
        For each $(q-2)$-photon $E$ in $\mathbb{R}^{q,q}$, there exists a unique pair of $(q-1)$-photons $(F,F')\in \SO_0(q,q)/P_q \times \SO_0(q,q)/P_q'$ such that $E=F\cap F'$.
    \end{prop}

    \begin{proof}
        Let $E$ be an isotropic subspace in $\R^{q,q}$ of dimension $q-1$. The orthogonal of $E$ is of dimension $q+1$ and is the orthogonal sum of $E$ with a subspace of signature $(1,1)$, i.e $E^{\perp} = E \oplus^{\perp} \R^{1,1}$. Taking an isotropic space of dimension $q$ containing $E$ is equivalent to taking an isotropic line in $\R^{1,1}$. There are exactly two isotropic lines in $\R^{1,1}$, hence the result.
    \end{proof}
    
    The space of $(q-2)$-photons is then the $(q',q)$-flag variety of $\SO_0(q,q)$. Every $(q-2)$-photon can be included in exactly two $(q-1)$-photons, one in each of the two flag varieties $\mathcal{F}^q = \SO_0(q,q)/P_q$ and $\mathcal{F}^{q'} = \SO_0(q,q)/P_{q'}$. This describes the projections onto the flag varieties associated to each of the last two roots. In order to make the notations more intuitive, we will often replace the pair $(q',q)$ in the notations with $q-1$. For instance we will write $\mathcal{F}^{q-1}$ instead of $\mathcal{F}^{q',q}$, or $P_{q-1}$ instead of $P_{q',q}$. With this notation the Borel subgroup of $\SO_0(q,q)$ is then $P_{1,\dots,q-1}$. We will also replace the data of a pair of $(q-1)$-photons in a flag variety containing $q'$ and $q$ by a single $(q-2)$-photon, as this data is equivalent according to Proposition \ref{prop:sub-max_photon}.

    Another property of flag varieties that is different from the case $p\neq q$ is that the flag varieties $\SO_0(q,q)/P_q$ and $\SO_0(q,q)/P_{q'}$ may not be self-opposite.

    \begin{prop}\label{maxphotonqq}
        The parabolic subgroups $P_q$ and $P_{q'}$ are self-opposite if and only if $q = 0$ mod $2$.
    \end{prop}
    
    \begin{corollary}
        When $q = 0$ mod $2$, every parabolic subgroup of $\SO_0(q,q)$ is self-opposite. When $q = 1$ mod $2$, a parabolic subgroup $P_\Theta $ of $\SO_0(q,q)$ is self-opposite if and only if $\Theta $ contains either neither or both the roots $q'$ and $q$.
    \end{corollary}
    
    \subsubsection{Parametrization of the flag varieties}
    
    The subgroup $P_\Theta $ (resp. $P_\Theta ^{opp}$) is the intersection of $\SO_0(p,q)$ with the set of block-lower triangular (resp. upper triangular) matrices with the square blocks on the diagonal of size $(i_1,i_2-i_1,\dots,i_k-i_{k-1}, p-q +2(q-i_k), i_k-i_{k-1},\dots, i_2-i_1)$. The intersection $L_\Theta =P_\Theta \cap P_\Theta ^{opp}$ is the subgroup of block-diagonal matrices with block sizes $(i_1,i_2-i_1,\dots,i_k-i_{k-1}, p-q +2(q-i_k), i_k-i_{k-1},\dots, i_2-i_1)$ and is called the Levi subgroup of $P_\Theta $. Let $U_\Theta $ be the unipotent radical of $P_\Theta ^{opp}$ and let $\mathfrak{u}_\Theta $ be the Lie algebra of $U_\Theta $. When $\Theta =(1,2,\dots,k-1,k)$, we will write an element of $\mathfrak{u}_\Theta $ in the following way:
    \renewcommand{\arraystretch}{1.3}
    \begin{equation}
    \left(
 			\begin{array}{ccccc|c|ccccc}
 				0 & -a_1^{k-1} & -a_1^{k-2} & \cdots & -a_1^1     & \bar v_1^0     & -b_1^1     & \cdots    & -b_1^{k-2} & -b_1^{k-1} & 0        \\
                ~ & 0          & -a_2^{k-2} & \cdots & -a_2^1     & \bar v_2^0     & -b_2^1     & \cdots    & -b_1^{k-1} & 0          & b_1^{k-1}\\
                ~ & ~          & 0          & \ddots & \vdots     & \vdots         & \vdots     & \iddots   & 0          & b_2^{k-2}  & b_1^{k-2}\\
                ~ & ~          & ~          & \ddots & -a_{k-1}^1 & \bar v_{k-1}^0 & -b_{k-1}^1 & \iddots   & \iddots    & \vdots     & \vdots   \\
                ~ & ~          & ~          & ~      & 0          & \bar v_k^0     & 0          & b_{k-1}^1 & \cdots     & b_2^1      & b_1^1    \\\hline
                ~ & ~          & ~          & ~      & ~          & 0_{p+q-2k}     & v_k^0      & v_{k-1}^0 & \cdots     & v_2^0      & v_1^0    \\\hline
                ~ & ~          & ~          & ~      & ~          & ~              & 0          & a_{k-1}^1 & \cdots     & a_2^1      & a_1^1    \\
                ~ & ~          & ~          & ~      & ~          & ~              & ~          & \ddots    & \ddots     & \vdots     & \vdots   \\
                ~ & ~          & ~          & ~      & ~          & ~              & ~          & ~         & 0          & a_2^{k-2}  & a_1^{k-2}\\
                ~ & ~          & ~          & ~      & ~          & ~              & ~          & ~         & ~          & 0          & a_1^{k-1}\\
                ~ & ~          & ~          & ~      & ~          & ~              & ~          & ~         & ~          & ~          & 0        \\
 			\end{array}\right)\label{eq:param}
\end{equation}

 		\renewcommand{\arraystretch}{1}
    where all $(a_j^i)_{1\leq j\leq k,~1\leq i\leq k-j}$ and $(b_j^i)_{1\leq j\leq k,~1\leq i\leq k-j}$ are real numbers, and the $(v_j^0)_{1\leq j\leq k}$ are (column) vectors in $\RR^{p-k,q-k}$. The line vectors $\bar v_j^0$ are defined such that $\bar v_j^0 v_j^0 = - Q(v_j^0)$. For any other choice of $\Theta =(i_1,\dots,i_k)$, $\mathfrak{u}_\Theta $ is a subspace of $\mathfrak{u}_{1,\dots,i_k}$ where a subset of the $a_j^i$ vanish. For any $1\leq j\leq k$ and for any $0\leq i\leq k-j$ we define $$v^i_j = (b^{i}_j,b^{i-1}_j,\dots, b^{1}_j,v_j^0,a^{1}_j,\dots,a^{i-1}_j,a^{i}_j)\in \RR^{p-k+i-1,q-k+i-1}.$$ We will also sometime denote the first (resp. the last) coordinate of $v_j^0$ by $b_j^0$ (resp. $a_j^0$), and the vector obtained by removing the fist and the last entries of $v_j^0$ by $v_j^{-1}$. For $1\leq j\leq k$, let $u_j(v_j^{k-j})$ be the square matrix of size $p+q$ whose $j$-th line is $(0,\dots,0,-a^{k-j}_j,\dots, -a^{1}_j,\bar v_j^0,-b^{1}_j,\dots,-b^{k-j}_j,0,\dots,0)$ and whose $(p+q+1-j)$-th column is $(0,\dots,0,b^{k-j}_j,\dots, b^{1}_j,v_j^0,a^{1}_j,\dots,a^{k-j}_j,0,\dots,0)$, with all other entries being $0$. A matrix in $\mathfrak{u}_\Theta $ then uniquely writes as $$u = u_1(v_1^{k-1}) +u_2(v_2^{k-2})+\dots + u_k(v_k^{0}).$$

    \begin{proposition}\label{prop:matrixproduct}
        The map $$\Psi_\Theta  : \begin{array}{rcl}
            \mathfrak{u}_\Theta  & \to & U_\Theta  \\
             u = u_1(v_1^{k-1})+\dots + u_k(v_k^{0}) &\mapsto & U= \exp\left(u_1(v_1^{k-1})\right)\dots\exp\left(u_k(v_k^{0})\right)
        \end{array}$$
        is a diffeomorphism.
    \end{proposition}

    \begin{proof}
        First note that for all $1\leq j\leq k$, the subspace
        $$V_j = \left\lbrace u_i(v_j^{k-j})~|~ v_j^{k-j}\in\RR^{p-j,q-j}\right\rbrace$$
        is an abelian Lie subalgebra of $\mathfrak{u}_\Theta $. Let $f_j =(f_j^1,\dots, f_j^{p+q-2j})$ be a basis of $V_j$. Then if $v_j^{k-j} = \sum \lambda_i f_j^i$, we have
        $$\exp\left(u_j(v_j^{k-j})\right) = \exp(\lambda_1 f_j^1))\dots \exp (\lambda_{p+q-2j} f_i^{p+q-2j}).$$
        Let $f = f_1\cup\dots\cup f_k$ be the corresponding basis of $\mathfrak{u}_\Theta $ as a vector space. This basis has the property that for any $1\leq i\leq \dim \mathfrak{u}_k$, the first $i$ vectors of $f$ span a subalgebra of $\mathfrak{u}_k$. Then by \cite{corwin1990representations} Prop. 1.2.8, the map $\Psi_\Theta $ is a diffeomorphism.
    \end{proof}

    \subsection{Transversality}\label{sec:transversality}
    
    \begin{definition}
        Let $\Theta  = (i_1,\dots,i_k)$ such that $\SO_0(p,q)/P_{\Theta}$ is self-opposite. Two flags $F_1,F_2\in\mathcal{F}^\Theta $ are called \emph{transverse} if for all $1\leq \ell\leq k$, $F_1^{i_\ell} \oplus (F_2^{i_\ell})^\perp = \RR^{p,q}$.
    \end{definition}

    The two standard flags $F_0$ and $F_\infty$ are transverse. Given a flag $F\in \mathcal{F}^\Theta $, we will denote by $\Omega(F)\subset \mathcal{F}^\Theta $ the subset of flags that are transverse to $F$.

    \begin{proposition}\label{prop:affine_chart}
        Let $\Theta  = (i_1,\dots,i_k)$ such that $\SO_0(p,q)/P_{\Theta}$ is self-opposite. For any flag $F\in\Omega(F_\infty^\Theta )$ transverse to $F_\infty^\Theta $ there exists a unique $U\in U_\Theta $ such that $F = U\cdot F^\Theta _0$. Thus $\Omega(F_\infty^\Theta )$ is diffeomorphic to $U_\Theta $, hence diffeomorphic to $\mathfrak{u}_\Theta $. In other words $\Omega(F_\infty^\Theta )$ is an affine chart of $\mathcal{F}^\Theta $, with a natural choice of origin given by $F^\Theta _0$.
    \end{proposition}

    Let $U\in U_\Theta $. The flag $F = U\cdot F_0$ is transverse to $F_\infty$ and we now want to express under which conditions it is transverse to $F_0$. Since for any $i\in\Theta $ we have $F_0^i = \Span(\Tilde{e}_1,\dots,\Tilde{e}_i)$, its orthogonal is $$\left(F_0^i\right)^\perp = \Span(\Tilde{e}_1,\dots,\Tilde{e}_q,x_{q+1},\dots,x_{p-q}, e_q,\dots, e_{i-1}).$$
    For $F^i = U\cdot F_0^i$ to be transverse to $F_0^i$, the projection of $F^i$ to $\Span(e_1,\dots,e_i)$ parallel to $\left(F_0^i\right)^\perp$ must be surjective. This corresponds to the matrix $U\in U_\Theta $ having its $i\times i$ upper right minor to be non-zero. We will denote this minor by $\det_i(U)$ in the following.
    We thus obtain the first description of $\Omega(F^\Theta _0)\cap \Omega(F^\Theta _\infty)$:
    $$\Omega(F^\Theta _0)\cap \Omega(F^\Theta _\infty) \overset{homeo}{\simeq} \left\lbrace U\in U_\Theta ~ |~ \forall i\in\Theta ,~\det_i(U)\neq 0\right\rbrace.$$
    
    \begin{rem}
        The space of triples of transverse flags is a fiber bundle over the space of pairs of transverse flags $(F_0, F_{\infty})$ with fiber $\Omega(F_0) \cap \Omega(F_{\infty})$. While we will focus on counting the connected components in $\Omega(F_0) \cap \Omega(F_{\infty})$ for each pair, since the connected group $\SO_0(p,q)$ acts transitively on the space of pairs of transverse flags, the space of triples of transverse flags has the same number of connected components as the quotient of $\Omega(F_0) \cap \Omega(F_{\infty})$ by the stabilizer $P_\Theta\cap P_\Theta^{opp}$ of the pair $(F_0,F_\infty)$.
    \end{rem}

    \subsection{Einstein universe}\label{sec:einstein}

    Let $p\geq q\geq 2$. The space $\SO_0(p,q)/P_1$ of isotropic lines in $\RR^{p,q}$ is called \emph{Einstein universe} and will be denoted by $\Ein_{p-1,q-1}$. This is a smooth manifold of dimension $p+q-2$ embedded in $\mathbb RP^{p+q-1}$, on which the quadratic form $Q$ induces a conformal class of pseudo-Riemannian metric of signature $(p-1,q-1)$. When $q=1$ the Einstein universe $\Ein_{p-1,0}$ is identified the boundary of the hyperbolic space $\mathbb H^{p}$ which is a sphere $\mathbb S^{p-1}$ with its conformal structure. When $q=2$, the space $\Ein_{p-1,1}$ is called \emph{Lorentzian}. Since Einstein universe only has a conformal structure, only the sign of the quadratic form on a tangent vector is well defined. In particular, the type (time, space or light) of a tangent vector is well defined, and so is the lightcone of a point. The lightcone of a point $x\in\Ein_{p-1,q-1}$ is given by:
    $$C(x) = \mathbb{P}_{\Iso }(x^\perp) = \mathbb{P}\lbrace v\in x^\perp~|~Q(v)=0\rbrace.$$
    The complement in $\Ein_{p-1,q-1}$ of any lightcone is an affine chart called \emph{Minkovski space}, denoted by $\Mink_{p-1,q-1}$. A choice of any point in $\Mink_{p-1,q-1} = \Ein_{p-1,q-1}\backslash C(x)$ identifies it with $\mathbb R^{p-1,q-1}$.  
    \begin{proposition}\label{prop:intersection_lightcones}
        Let $q\geq 2$. Let $x,x'\in\Ein_{p-1,q-1}$, with $x'\notin C(x)$. The space of photons passing through $x$ is identified with $C(x)\cap C(x')$ which is conformally equivalent to $\Ein_{p-2,q-2}$.
    \end{proposition}

    \begin{proof}
        The set of photons passing through $x$ is the set of isotropic planes in $\RR^{p,q}$ containing the line $x$. Let $\Delta$ be a photon passing through $x$. Since $x\notin C(x')$, $\Delta$ is not contained in $x'^\perp$. However since $\dim x'^\perp = p+q-1$, the intersection of $\Delta$ and $x'^\perp$ is an isotropic line in $x^\perp\cap x'^\perp$, hence defines a point in $C(x)\cap C(x')$. Conversely, any point $y\in C(x)\cap C(x')$ gives rise to a unique photon passing through $x$ defined by $\Delta = x\oplus y$. Since the quadratic form $Q$ restricted to the space $x^\perp\cap x'^\perp$ is of signature $(p-1,q-1)$, the set $C(x)\cap C(x')$ of isotropic lines in $x^\perp\cap x'^\perp$ is identified with $\Ein_{p-2,q-2}$.
    \end{proof}
    
    \begin{prop}
        The space $\widehat{\Ein}_{p,q}$ of isotropic half-lines is a double cover of $\Ein_{p,q}$ which is conformally equivalent to the pseudo-Riemannian space $(\mathbb{S}^p \times \mathbb{S}^q, ds_p^2 - ds_q^2)$. In such a conformal model, the $k$-dimensional photons of $\widehat{\Ein}_{p,q}$ are exactly the graphs of isometries $f : \mathbb{S}^k \subset \mathbb{S}^p \rightarrow \mathbb{S}^q$ where $\mathbb{S}^k \subset \mathbb{S}^p$ is a totally geodesic embedding.
    \end{prop}
    
    \begin{proof}
        Let $\R^{p+1, q+1} = \R^{p+1,0} \oplus \R^{0, q+1}$ be an orthogonal splitting and let $S^+$ and $S^-$ be the unit spheres of $\R^{p+1,0}$ and $\R^{0, q+1}$. Finally, let us consider the map
        
        \[\varphi : (x^+, x^-) \in S^+ \times S^- \longmapsto [x^+ + x^-] \in \widehat{\Ein}_{p,q}.\]
        
        The map $\varphi$ is a diffeomorphism since for each half-line in $\widehat{\Ein}_{p,q}$ there exists unique points $(x^+, x^-) \in S^+ \times S^-$ such that $[x^+ + x^-]$ is this half-line. When one endows $S^+ \times S^-$ with the pseudo-Riemannian metric $ds_p^2 - ds_q^2$ where $ds_p^2$ (resp. $ds_q^2$) is the round Riemannian metric induced by $\R^{p+1,0}$ (resp. $\R^{0, q+1}$), the map $\varphi$ becomes a conformal equivalence as the image of the map $(x^+, x^-) \mapsto x^+ + x^-$ in the isotropic cone of $\R^{p+1, q+1}$ is transverse to all isotropic half-lines. 
        
        We then only have to check that the $k$-dimensional photons in $\mathbb{S}^p \times \mathbb{S}^q$ are indeed the graphs of isometries $f : \mathbb{S}^k \subset \mathbb{S}^p \rightarrow \mathbb{S}^q$. Let $V$ be an isotropic subspace of dimension $k+1$. One may embed $V$ into $\R^{k+1, k+1}$ by adding a transverse isotropic space, thus we may assume that $p=q=k$. Since $V$ is in direct sum with both $\R^{k+1,0}$ and $\R^{0, k+1}$, it can be seen as the graph of a linear map $f : \R^{k+1,0} \rightarrow \R^{0, k+1}$. This graph is isotropic if and only if $f$ is an isometry between $\R^{k+1,0}$ and $- \R^{0, k+1}$, which yields the intended result.
    \end{proof}

        The flag variety $\SO_0(q,q)/P_q$ corresponds to the orientation preserving isometries of $\mathbb{S}^{q-1}$ and $\SO_0(q,q)/P_{q'}$ to orientation reversing isometries. With this geometric description we can see why the parabolic subgroups $P_{q'}$ and $P_q$ of $\SO_0(q,q)$ may only be self-opposite when $q$ is even. Let $F_0 \in \SO_0(q,q)/P_q$ be the graph of the identity map in $\mathbb{S}^{q-1}$ for a conformal model $\widehat{\Ein}_{q-1, q-1} \simeq \mathbb{S}^{q-1} \times \mathbb{S}^{q-1}$. Two maximal photons in $\widehat{\Ein}_{q-1, q-1}$ are transverse if and only if they do not intersect. Any other maximal photon in $\widehat{\Ein}_{q-1, q-1}$ can be seen as the graph of an isometry $g : \mathbb{S}^{q-1} \rightarrow \mathbb{S}^{q-1}$. When $q$ is odd, any direct isometry in $\SO_0(q)$ must have a fixed point, meaning that its graph in $\widehat{\Ein}_{q-1, q-1}$ must intersect $F_0$. Therefore any maximal photon transverse to $F_0$ must be within $\SO_0(q,q)/P_{q'}$. Conversely, when $q$ is even, any maximal photon not intersecting $F_0$ must be within $\SO_0(q,q)/P_q$. This yields the following result~:
        
        \begin{prop}
            The space of $q$-flags for $\SO_0(q,q)$ is the union of two connected components corresponding to the flag varieties $\SO_0(q,q)/P_q$ and $\SO_0(q,q)/P_{q'}$. Those two flag varieties are self-opposite if and only if $q$ is even.
        \end{prop}
        
    In this case, the two spaces $\SO_0(q,q)/P_q$ and $\SO_0(q,q)/P_{q'}$ are isomorphic to each other~: 
    
    \begin{prop}\label{transpreserving}
        When $q$ is even, there is a transversality-preserving diffeomorphism between $\SO_0(q,q)/P_q$ and $\SO_0(q,q)/P_{q'}$.
    \end{prop}
    
    \begin{proof}
        Let $s$ be an orientation reversing isometry and let us consider the map $\varphi : f \in \SO(q) \mapsto s \circ f$. It is a diffeomorphism between $\SO_0(q,q)/P_q$ and $\SO_0(q,q)/P_{q'}$ as it sends orientation preserving isometries to orientation reversing ones. For two maps $f,g \in \SO(q+1)$, there exists a point $x$ in $\mathbb{S}^q$ such that $f(x) = g(x)$ if and only if there exists a point $x$ in $\mathbb{S}^q$ such that $s \circ f(x) = s \circ g(x)$, thus the map $\varphi$ preserves transversality.
    \end{proof}

    \section{Equations of transversality}
    
    The goal of this section is to provide a description of the space $\Omega(F_0^\Theta )\cap\Omega(F_\infty^\Theta )$ of $\Theta $-flags transverse to both $F_0^\Theta $ and $F_\infty^\Theta $. We gave a first description of that space at the end of Section \ref{sec:transversality} as the subset of matrices in $U_\Theta $ having all their top-right $i\times i$ minors to be non-zero, for all $i\in\Theta $. We now want to give an explicit computation of these minors, in terms of the parametrization \ref{eq:param} of $U_\Theta $. This involves computing the image $U$ of an element in $\mathfrak{u}_\Theta $ by the map $\Psi_\Theta $, and then computing the minors $\det_i(U)$. While we carry the general computation in Section \ref{sec:general_eq_transversality}, we first provide detailed computations for the two most elementary non-trivial cases, namely $(1,2)$-flags of $\SO_0(p,q)$ (Section \ref{sec:transversality_pointed_photons}) and $(1,2,3)$-flags in $\SO_0(p,q)$ (Section \ref{sec:123-photon_eq_transversality}). We start with these two special cases to demonstrate how the methods we use to carry these computation correspond to a geometric interpretation of relative positions of photons in the Einstein universe.

    \subsection{\texorpdfstring{$(1,2)$}{(1,2)}-flags}\label{sec:transversality_pointed_photons}

    Let $p\geq q\geq 2$. We will consider the space $\mathcal{F} = \mathcal{F}^{1,2}$ of $(1,2)$-flags (also called pointed photons) in $\RR^{p,q}$. Let $F_0 = (x_0,\Delta_0)$ and $F_\infty = (x_\infty,\Delta_\infty)$ be the two standard flags defined in Section \ref{sec:transversality}. From Proposition \ref{prop:affine_chart}, any pointed photon $F =(x,\Delta)\in\mathcal{F}$ transverse to $F_\infty$ can be written as $U\cdot(x_0,\Delta_0)$ where 
    \begin{align*}
        U &= \exp\left(\begin{array}{ccccc}
         0 &-a & \bar v_1 & -b & 0\\
         ~ & 0 & 0 & 0  & b\\
         ~ & ~ & 0  & 0 & v_1\\
         ~ & ~ & ~          & 0 & a\\
         ~ & ~ & ~          & ~ & 0
    \end{array}\right)\exp\left(\begin{array}{ccccc}
         0 & 0 & 0 & 0 & 0\\
         ~ & 0 & \bar v_2 & 0  & 0\\
         ~ & ~ & 0  & v_2 & 0\\
         ~ & ~ & ~          & 0 & 0\\
         ~ & ~ & ~          & ~ & 0
    \end{array}\right)\\
    &= \left(\begin{array}{ccccc}
         1 &-a & \bar v_1 - a\bar v_2 & -b +\frac{a}{2}Q(v_2)-B(v_1,v_2) & -\frac{Q(v_1^1)}{2}\\
         ~ & 1& \bar v_2 & -\frac{Q(v_2)}{2}  & b\\
         ~ & ~ & I & v_2 & v_1\\
         ~ & ~ & ~       & 1 & a\\
         ~ & ~ & ~       & ~ & 1
    \end{array}\right).
    \end{align*}
    
    Compared to the parametrization \ref{eq:param}, to alleviate the notations we write in this section $v_1\in \RR^{p-2,q-2}$ (resp. $v_2$) instead of $v_1^0$ (resp. $v_2^0$), $a$ (resp. $b$) instead of $a_1^1$ (resp. $b_1^1$) and we have then $v_1^1 = (b,v_1,a)\in\RR^{p-1,q-1}$. For the flag $F$ to be transverse to $F_0$, we first need to have $x\oplus x_0^\perp = \RR^{p,q}$. In other words, $x$ needs not to lie on the lightcone of $x_0$ in $\Ein_{p-1,q-1}$. Since the coordinates of $x$ in the affine chart defined by $x_\infty$ and $x_0$ are precisely given by $v_1^1$, the transversality condition is $Q(v_1^1)\neq 0$ which is equivalent to $\det_1(U) = -\frac{Q(v_1^1)}{2}\neq 0$. 
    
    The second transversality condition is given by $\det_2(U)\neq 0$. We can interpret this condition geometrically as follows: there is a unique photon $\Delta'_0$ passing through $x$ and intersecting $\Delta_0$ since $\Delta_0$ intersect $C(x)$ in one point, which is given by $\Delta_0'=\Span(b\Tilde{e}_1 +Q(v_1^1)\Tilde{e}_2)$. Both $\Delta$ and $\Delta'_0$ are photons passing through $x$, so we can see them as two points in $E = C(x)\cap C(x_\infty)\simeq\Ein_{p-2,q-2}$ (see Proposition \ref{prop:intersection_lightcones}). The photons $\Delta$ and $\Delta_0$ are transverse if and only if $\Delta\cap C(x_\infty)$ does not lie on the lightcone of $\Delta'_0\cap C(x_\infty)$ in $E$. Note that $E$ is endowed with a natural affine chart where the lightcone of $\Delta_\infty$ is at infinity. The point $\Delta'_0\cap C(x_\infty)\in E$ is on the lightcone of $\Delta_\infty\cap E$ precisely when $Q(v_1)=0$. Otherwise, it lies in the affine chart $E\backslash C(\Delta_\infty\cap E)$. When $Q(v_1)\neq 0$, we use the affine chart $\RR^{p-2,q-2} = E\backslash C(\Delta_\infty\cap E)$ with origin $\Delta'_0\cap E$. In this chart, the coordinate of $\Delta\cap E$ is $v_2^{(1)} = v_2+\frac{2b}{Q(v_1)}v_1$. This means that when $Q(v_1)\neq 0$, the transversality condition on $\Delta$ is $Q(v_2^{(1)}) \neq 0$. We now look at the equation given by $\det_2(U)\neq 0$:
    $$\begin{array}{rcl}
        \mathrm{det}_2(U) & = &\det\left(\begin{array}{cc}
        -b +\frac{a}{2}Q(v_2)-B(v_1,v_2) & -\frac{Q(v_1^1)}{2} \\
        -\frac{Q(v_2)}{2}  & b
    \end{array}\right) \\&\overset{L_1\leftarrow L_1+aL_2}{=} &\det\left(\begin{array}{cc}
        -b -B(v_1,v_2) & -\frac{Q(v_1)}{2} \\
        -\frac{Q(v_2)}{2}  & b
    \end{array}\right) \\& \overset{L_1\leftarrow L_2+\frac{2b}{Q(v_1)}L_1}{=} &\det\left(\begin{array}{cc}
        -b -B(v_1,v_2) & -\frac{Q(v_1)}{2} \\
        -\frac{1}{2} Q(v_2+\frac{2b}{Q(v_1)}v_1) & 0
    \end{array}\right)
    \end{array}$$

    The equation $\det_2(U)\neq 0$ is then $\frac{1}{4}Q(v_1)Q(v_2+\frac{2b}{Q(v_1)}v_1)\neq 0$ and we retrieve the equation obtained via geometric arguments. Note however that $Q(v_1)$ may vanish while $(x,\Delta)$ is still transverse to $(x_0,\Delta_0)$. In that case, the point $\Delta'_0\cap E$ does not lie in the affine chart $E\backslash C(\Delta_\infty\cap E)$. The intersection of the lightcone of $\Delta'_0\cap E$ with the affine chart is then a hyperplane for which we will not need to give an explicit euqation.

    \subsection{\texorpdfstring{$(1,2,3)$}{(1,2,3)}-flags}\label{sec:123-photon_eq_transversality}
    
    Let $p \geq q\geq 3$. We will now consider the space $\mathcal{F} = \mathcal{F}^{1,2,3}$ of $(1,2,3)$-flags. Let $F_0 = (x_0,\Delta_0,\Phi_0)$ and $F_\infty = (x_\infty,\Delta_\infty,\Phi_\infty)$ be the two standard flags defined in Section \ref{sec:transversality}. From Proposition \ref{prop:affine_chart}, any $(1,2,3)$-flag $F =(x,\Delta,\Phi)\in\mathcal{F}$ transverse to $F_\infty$ can be written as $U\cdot(x_0,\Delta_0,\Phi_0)$ where 
    $$U = \exp\left(u_1(v_1^2)\right)\exp\left(u_2(v_2^1)\right)\exp\left(u_3(v_3^0)\right).$$
    Now we observe that the fibers of the map $$\begin{array}{rcl}
        \mathcal{F}^{1,2,3} & \to & \Ein_{p-1,q-1} \\
        (x,\Delta,\Phi) & \mapsto & x
    \end{array}$$
    above a point $x$ transverse to $x_0$ is identified to the space of pointed photons in $\Ein_{p-2,q-2}$. Our goal is to use the description given in the Section \ref{sec:transversality_pointed_photons} to describe this fiber.

    We start with $x = U\cdot x_0$ which is transverse to $x_0$ exactly when $Q(v_1^2)\neq 0$, as seen in Section \ref{sec:transversality_pointed_photons}. We now want to see $\Delta_0$ and $\Phi_0$ in the space $E = C(x)\cap C(x_\infty)\simeq \Ein_{p-2,q-2}$. As in Section \ref{sec:transversality_pointed_photons}, we define $\Delta'_0$ to be the only photon passing through $x$ that intersect $\Delta_0$, and in the same spirit we define $\Phi'_0$ to be the set of all photons passing through $x$ and intersecting $\Phi_0$. To alleviate the notations we will write $\Delta'_0$ (resp. $\Delta$, $\Phi$, $\Phi'_0$, $\Delta_\infty$, $\Phi_\infty$) instead of their respective intersection with $E$. The point $\Delta'_0$ lies on the lightcone of $\Delta_\infty$ exactly when $Q(v_1^1) =0$, and $\Phi'_0$ intersects $\Phi_\infty$ exactly when $Q(v_1^0)\neq 0$. In the following, we assume $Q(v_1^1)\neq 0$ and $Q(v_1^0) \neq 0$. We now have two pointed photons $(\Delta,\Phi)$ and $(\Delta_0', \Phi_0')$ in $E\simeq \Ein_{p-2,q-2}$, and we want to apply Section \ref{sec:transversality_pointed_photons} computations. For this we need to use a basis  $e'=(e'_2,\dots,e'_{q},x'_{q+1},\dots,x'_{p-q},\Tilde{e}'_{q},\dots,\Tilde{e}'_2)$ of $x^\perp \cap x_\infty^\perp= (e_2,\dots, \Tilde{e}_2)\subset \RR^{p,q}$ in which $\Delta'_0 = \Span(\Tilde{e}'_2)$, $\Phi'_0 =\Span(\Tilde{e}'_2,\Tilde{e}'_3)$, $\Delta_\infty = \Span(e'_2)$ and $\Phi_\infty = \Span(e'_2,e'_3)$ (i.e. in which $(\Delta_0', \Phi_0')$ and $(\Delta_\infty,\Phi_\infty)$ are the standard $(1,2)$-flags). Such a basis is obtained by applying the matrix
    \renewcommand{\arraystretch}{2}
    $$ \left(\begin{array}{ccccc}
         1 &-\dfrac{2b_1^2a_1^{1}}{Q(v_1^1)} & \dfrac{-2b_1^2}{Q(v_1^1)}\bar v_1^0 & -\dfrac{2b_1^2b_1^1}{Q(v_1^0)}   & -\dfrac{4(b_1^2)^2}{Q(v_1^0)}\\
         ~ & 1                               & \dfrac{-2b_1^1}{Q(v_1^0)}\bar v_1^0 & -\dfrac{4(b_1^1)^2}{Q(v_1^0)}    & \dfrac{2b_1^1b_1^2}{Q(v_1^1)}\\
         ~ & ~                               & I                                   & \dfrac{-2b_1^1}{Q(v_1^0)}v_1^0   & \dfrac{-2b_1^2}{Q(v_1^1)}v_1^0\\
         ~ & ~                               & ~                                   & 1                                & \dfrac{2b_1^2a_1^{1}}{Q(v_1^1)}\\
         ~ & ~                               & ~                                   & ~                                & 1
    \end{array}\right) \in \SO_0(p-1,q-1)$$
        \renewcommand{\arraystretch}{1}
    to the basis $(e_2,\dots, \Tilde{e}_2)$.
    The coordinates of $\Delta$ in the affine chart $E\backslash C(\Delta_\infty)$ are then
    $$\Delta = b_2^{1,(1)}e'_3 + v_2^{0,(1)} + a_2^{1,(1)}\Tilde{e}'_3,$$
    where 
    $$ a_2^{1,(1)} = a_2^1 - \frac{2b_1^2}{Q(v_1^0)}a_1^1, $$
    $$ v_2^{0,(1)} = v_2^0 - \frac{2b_1^2}{Q(v_1^0)}v_1^0, $$
    $$ b_2^{1,(1)} = b_2^1 - \frac{2b_1^2}{Q(v_1^0)}b_1^1 - \frac{2(b_1^1)^2a_2^1}{Q(v_1^0)} - \frac{2b_1^1}{Q(v_1^0)}B(v_1^0,v_2^0). $$
    
    and $\Phi = \Span(\Delta, -\frac{Q(v_3^{0,(1)})}{2}e'_3 + v_3^{0,(1)} + \Tilde{e}'_3)$ where 
    $$ v_3^{0,(1)} = v_3^0 - \frac{2b_1^2}{Q(v_0^1)}v_1^0 .$$
    
    We can now apply the computation done in Section \ref{sec:transversality_pointed_photons} to the triple of $(1,2)$-flags $(\Delta_0', \Phi_0')$, $(\Delta,\Phi)$ and $(\Delta_\infty,\Phi_\infty)$. That triple will be transverse when $v_3^{0,(2)} = v_3^{0,(1)} +\frac{2b_2^{1,(1)}}{Q(v_2^{0,(1)})}v_2^{0,(1)}$ has non-zero norm.

    \subsection{General case}\label{sec:general_eq_transversality}

    Let $p \geq q \geq 2$ and $1 \leq k \leq q$. Let $\mathcal{F}=\mathcal{F}^{1,\dots,k}$. For $F^k_0$ and $F^k_{\infty}$ the two standard $(k-1)$-photons, the equation of transversality of $F = U \cdot F_0$ in the affine chart defined by $F_{\infty}$ and $F_0$ is $\det_k(U) = (-1)^k\det(S) \neq 0$ where $S = (s_{i,j})$ is the anti-transpose (i.e. the transpose with respect to the anti-diagonal) of the upper-right $k\times k$ submatrix of $U$. Computing the image of the map $\psi$ defined in \ref{prop:matrixproduct}, we obtain $s_{i,j} = b_j^{k-i+1}$ when $i > j$, $s_{i,i} = - \frac{1}{2} Q(v_i^{k-i})$ and 

    \[s_{i,j} = - b_i^{k-j+1} - Q(v_i^{k-j+1}, v_j^{k-j+1}) - \sum_{l=1}^{j-i} a_i^{k-l} s_{i+l, j}\]
    when $i < j$.
    We are going to perform a succession of change of variables. In order to keep coherent notations, the initial parameters described in \ref{eq:param} will be denoted with an exponent $(0)$. Let us perform elementary operations on the matrix $S$ in order to eliminate the coefficients on the first column starting from the second row, assuming that for all $i$, $Q(v_1^{i, (0)}) \neq 0$. This is just a variation of the standard Gaussian elimination, in which we weave in transformations of the first line so that the computations are easier. We apply the following transformations~:

    \[ \begin{split} L_1 &\longleftarrow L_1 + a_1^{k-1, (0)}L_2,  \\
    L_2 &\longleftarrow L_2 + \frac{2 b_1^{k-1, (0)}}{Q(v_1^{k-2, (0)})} L_1, \\
    &... \\
    L_1 &\longleftarrow L_1 + a_1^{1, (0)} L_k, \\
    L_k &\longleftarrow L_k + \frac{2 b_1^{1, (0)}}{Q(v_1^{0, (0)})} L_1.
    \end{split}\]

    We then obtain the matrix $S'$ where for all $i > 1$, $s'_{i,1} = 0$. Our goal will be to introduce new variables $(a_j^{i,(1)}),(b_j^{i,(1)}), (v_j^{i,(1)})$ such that the sub-matrix $S_2$ defined by $S_2 = (S')_{i >1, j > 1}$ admits the same expression as the upper-left $(k-1)\times (k-1)$ submatrix of $S$ after the change of variable $v_j^{i, (0)} \leftarrow v_j^{i, (1)}$, $a_j^{i, (0)} \leftarrow a_j^{i, (1)}$ and $b_j^{i, (0)} \leftarrow b_j^{i, (1)}$ so that we can proceed by induction.
    
    For all $i > j$, let us define 
    
    \[ \begin{split} b_j^{k-i+1, (1)} &= s'_{i,j} \\
    &= b_j^{k-i+1, (0)} - \frac{2 b_1^{k-i+1, (0)}}{Q(v_1^{k-i+1, (0)})} b_1^{k-j+1, (0)} - \frac{2 b_1^{k-i+1, (0)}}{Q(v_1^{k-i+1, (0)})} Q(v_1^{k-i-1, (0)}, v_j^{k-i-1, (0)})  \\
    &- \sum_{l=1}^{i-j} \frac{2 b_1^{k-i+1, (0)} a_j^{k-j-l, (0)} b_1^{k-j-l, (0)}}{Q(v_1^{k-i+1, (0)})} \end{split}\]
    
    and $v_i^{k-i, (1)} = v_i^{k-i, (0)} + \frac{2 b_1^{k-i+1, (0)}}{Q(v_1^{k-i, (0)})} v_1^{k-i, (0)}$. We can then see that

    \[s'_{i,i} = - \frac{1}{2} Q(v_i^{k-i, (0)}) - \frac{2 b_1^{k-i+1, (0)}}{Q(v_1^{k-i, (0)})} - \frac{2 (b_1^{k-i+1, (0)})^2}{Q(v_1^{k-i, (0)})} = - \frac{1}{2} Q(v_i^{k-i, (1)}).\]

    Let us then define $a_j^{k-j-1, (1)} = a_j^{k-j-1, (0)} + \frac{2 b_1^{k-1, (0)} a_j^{k-j-1, (0)}}{Q(v_1^{k-2, (0)})}$ and for $i < k-j-1$, 

    \[a_j^{i, (1)} = a_j^{i, (0)} + \frac{2 b_1^{k-j+1, (0)} a_j^{i, (0)}}{Q(v_1^{k-j, (0)})} + \sum_{l=i+1}^{k-j} \frac{2 a_j^{l, (1)} b_1^{l, (0)} a_1^{i, (0)}}{Q(v_1^{l-1, (0)})},\]

    as well as, for $i \leq k-j$, 

    \[v_j^{i, (1)} = v_j^{i, (0)} + \frac{2 b_1^{k-j+1, (0)}}{Q(v_1^{k-j, (0)})} v_1^{i, (0)} + \sum_{l=i}^{k-j-1} \frac{a_j^{l, (1)} b_1^{l, (0)}}{Q(v_1^{l, (0)})} v_1^{i, (0)}.\]

    By some elementary but tedious computations which we will not explicit, we then get that for $j > i$, 

    \[s'_{i,j} = - b_i^{k-j+1, (1)} - Q(v_i^{k-j+1, (1)}, v_j^{k-j+1, (1)}) - \sum_{l=1}^{j-i} a_i^{k-l, (1)} s'_{i+l, j}\]

    and that for all $i \leq k-j$, 

    \[Q(v_j^{i, (1)}) = Q(v_j^{i-1, (1)}) + 2 a_j^{i,(1)} b_j^{i, (1)},\]

    thus showing that the submatrix $S_2$ is indeed of the same form as $S$ with the change of variable $v_j^{i, (0)} \leftarrow v_j^{i, (1)}$, $a_j^{i, (0)} \leftarrow a_j^{i, (1)}$ and $b_j^{i, (0)} \leftarrow b_j^{i, (1)}$. When for all $0\leq i\leq k-1$, $Q(v_1^{i, (0)})$ is assumed to be non-zero, the equation describing transversality to $F^k_0$ then writes $(-1)^k \frac{Q(v_1^{0, (0)})}{2}\det(S_2) \neq 0$ and we can re-apply this process, thus getting new variables $v_j^{i, (m)}$ for $1 \leq m \leq k-1$. Assuming that for all $i,j,m$, $Q(v_j^{i, (m)})$ is non-zero, we get the new variables for $i \leq k-j$ and $m \leq j-1$~:

    \[ \begin{split}
        v_j^{i, (m+1)} &= v_j^{i, (m)} + \frac{2 b_{m+1}^{k-m-j+1, (m)}}{Q(v_{m+1}^{k-m-j, (m)})} v_{m+1}^{i, (m)} + \sum_{l=i}^{k-m-j-1} \frac{a_j^{l, (m+1)} b_{m+1}^{l, (m)}}{Q(v_{m+1}^{l, (m)})} v_{m+1}^{i, (m)}, \\
        a_j^{i, (m+1)} &= a_j^{i, (m)} + \frac{2 b_{m+1}^{k-m-j+1, (m)} a_j^{i, (m)}}{Q(v_{m+1}^{k-m-j, (m)})} + \sum_{l=i+1}^{k-m-j} \frac{2 a_j^{l, (m+1)} b_{m+1}^{l, (m)} a_{m+1}^{i, (m)}}{Q(v_{m+1}^{l-1, (m)})}, \\
        b_j^{i, (m+1)} &= b_j^{i, (m)} - \frac{2 b_{m+1}^{i, (m)}}{Q(v_{m+1}^{i, (m)})} b_{m+1}^{k-j+1, (m)} - \frac{2 b_{m+1}^{i, (m)}}{Q(v_{m+1}^{i, (m)})} Q(v_{m+1}^{i, (m)}, v_j^{i, (m)}) \\
        &- \sum_{l=1}^{k-m-i+1-j} \frac{2 b_{m+1}^{i, (m)} a_j^{k-m-j-l, (m)} b_{m+1}^{k-m-j-l, (m)}}{Q(v_{m+1}^{i, (m)})}.
    \end{split}\]

    In the end, the equations of transversality for the $(1,\dots,k)$-flags become
    \[ \begin{split}
        \det_1(U) &= \frac{-1}{2}Q(v_1^{k-1, (0)}) \neq 0, \\
        \det_2(U) &=\frac{1}{4}Q(v_1^{k-2, (0)})Q(v_2^{k-2, (1)}) \neq 0, \\
        ...& \\
        \det_k(U) &=\frac{(-1)^{k+1}}{2^k}Q(v_1^{0, (0)})Q(v_2^{0, (1)})...Q(v_k^{0, (k-1)}) \neq 0. \\
    \end{split}\]

    \subsection{Realization as minors}

    Let $(F_0, U \cdot F_0, F_\infty)$ be a triple of transverse $(1,\dots,k)$-flags, for $k\leq q$. In order to determine the connected component of $\Omega(F_0)\cap\Omega(F_\infty)$ in which $U\cdot F_0$ lies, we will need the following data: 
    \begin{itemize}
        \item For all $1\leq j\leq k$ and $0\leq i\leq k-j$, we need to know the sign of $Q(v_j^{i,(j-1)})$. Moreover, if $k-i = q-1$ (this happens when $k = q-1$ and $i=0$ or when $k =q$ and $i=1$) or $k-i=q=p-1$ (this happens when $k=q=p-1$ and $i=0$) we need to know the sign of the first coordinate $(v_j^{i,(j-1)})_1$ of $v_j^{i,(j-1)}$ (to differentiate between future and past vectors when $k-i=q-1$ and between positive and negative numbers when $k-i=q=p-1$)
        \item For all $1\leq j\leq k-1$ and $1\leq i\leq k-j$, we need to know the sign of $b_j^{i,(j-1)}$.
    \end{itemize}

    The purpose of this section is to show that this data can be retrieved only using explicit minors of $U$. For $I,J$ two subsets of $\lbrace 1,\dots,p+q\rbrace$ of same cardinal, we will denote by $\Delta_{I,J}(U)$ the $(I,J)$-minor of $U$, i.e. the determinant of the submatrix $(u_{i,j})_{i\in I, j\in J}$. We then claim the following equalities: $\forall 1\leq j\leq k,~\forall 0\leq i\leq k-j,$
    $$\frac{(-1)^j}{2^j}Q(v_1^{i,(0)})Q(v_2^{i,(1)})\dots Q(v_j^{i,(j-1)}) = \Delta_{\lbrace 1,\dots,k-i \rbrace,\lbrace j+1,\dots,k-i,p+q-j+1,\dots,p+q\rbrace}(U)$$
    and $\forall 1\leq j\leq k,~\forall 1\leq i\leq k-j,$
    $$ \frac{(-1)^{j-1}}{2^{j-1}}Q(v_1^{i,(0)})Q(v_2^{i,(1)})\dots Q(v_{j-1}^{i,(j-2)})b_j^{i,(j-1)} = \Delta_{\lbrace 1,\dots, k-i+1\rbrace,\lbrace j,\dots,k-i ,p+q-j+1,\dots,p+q\rbrace}(U).$$
    In both cases, the quantities we are interested in are ratios of two minors of $U$. The data of all the signs of the numbers stated above is equivalent to the data of all the signs of the minors above.
    
    When $k-i = q-1$ or $k-i=q=p-1$, $\forall 1\leq j\leq k$:
    $$ \frac{(-1)^{j-1}}{2^{j-1}}Q(v_1^{i-1,(0)})\dots Q(v_{j-1}^{i-1,(j-2)})(v_j^{i-1,(j-1)})_1 = \Delta_{\lbrace 1,\dots, k-i+2\rbrace,\lbrace j,\dots,k-i+1 ,p+q-j+1,\dots,p+q\rbrace}(U)$$
    In both cases this is well defined, and $v_j^{-1,(j-1)}$ belongs to a space on which the restriction of $Q$ is positive definite, hence the sign of $(v_j^{i-1,(j-1)})_1$ is always the same as the one of $\Delta_{\lbrace 1,\dots, k-i+2\rbrace,\lbrace j,\dots,k-i+1 ,p+q-j+1,\dots,p+q\rbrace}(U)$. This observation is crucial because the number $$Q(v_1^{-1,(0)})Q(v_2^{-1,(1)})\dots Q(v_{j-1}^{-1,(j-2)})$$ does not arise as one of the previously computed minors.

    In the case when $k = q-1$ which correspond to the $\Theta$-positive structure of $\SO_0(p,q)$, the total number of minors to compute is $p(p-1)$, which is the same number as the $\Theta$-length of the longest word $w_\Theta^{max}$ of the $\Theta$-Weyl group defined in \cite{GW24}. This result is analogous to what happens in the case of full flags in a semi-simple split Lie group: the number of (generalized) minors required to determine the connected component is the length of the longest element of the Weyl group.

    \section{Counting connected components of triples of transverse flags}
    
    This section is dedicated to the count of connected components stated in Theorems \ref{th:CC_pq}, \ref{th:CC_qq+1} and \ref{th:CC_qq}. While we have a general method which would give the result for any flag variety of any group of the form $\SO_0(p,q)$ with $p\geq q$, in many cases much simpler arguments allow us to conclude. The most intricate case, in which we will need to apply this general method, is the case of $(1,\dots,q-1)$-flags in $\SO_0(p,q)$ with $p\neq q$ (see Section \ref{submaximalflags}). This general method works as follows: we use the parametrization \ref{eq:param}, and start with a subset of $\Omega(F_0)\cap\Omega(F_\infty)$ where all the $Q(v_j^{i,(j-1)})$ are non-zero. This subset is an open dense subset of $\Omega(F_0)\cap\Omega(F_\infty)$, for which counting the number of connected components--which we will call \emph{cells}--is straightforward (the exact number depend on the flag varieties and values of $p$ and $q$). 
    
    \begin{definition}
    A \emph{cell} is a connected subset of $\Omega(F_0)\cap\Omega(F_\infty)$ where all the $Q(v_j^{i,(j-1)})$ defined in \ref{sec:general_eq_transversality} are non-zero. The precise range of the indices $i$ and $j$ vary depending on which flag variety $F_0$ and $F_\infty$ belong to.
    \end{definition}
    
    \begin{remark}
        Note that since $Q(v_j^{i, (j-1)}) = Q(v_j^{i-1, (j-1)}) + 2 a_j^{i,(j-1)} b_j^{i, (j-1)}$, if a cell has opposite signs for $Q(v_j^{i, (j-1)})$ and $Q(v_j^{i-1, (j-1)})$, then the sign of $a_j^{i,(j-1)} b_j^{i, (j-1)}$ (hence of $b_j^{i, (j-1)}$) is constant inside it. A cell is then uniquely defined by the data of the signs of all the $Q(v_j^{i, (j-1)})$, plus the data of the signs of all $b_j^{i, (j-1)}$ for the $i,j$ such that $Q(v_j^{i, (j-1)})$ and $Q(v_j^{i-1, (j-1)})$ have opposite signs.
    \end{remark}
    
    These cells are homeomorphic to open balls of the same dimension as $\Omega(F_0)\cap\Omega(F_\infty)$. Our strategy revolves around identifying which of these cells lie in the same connected component of $\Omega(F_0)\cap\Omega(F_\infty)$ by identifying which one are separated by a codimension 1 locus of $\Omega(F_0)\cap\Omega(F_\infty)$ defined by $Q(v_j^{i,(j-1)})=0$ for given $i,j$. Once we have done this gluing operation for all $i,j$, we obtain the connected components of $\Omega(F_0)\cap\Omega(F_\infty)$, as well as the description of all the cells they contain.
    
    This section is divided in many different subsections, because of the many different behaviors of flag varieties depending on the values of $p$ and $q$, and on the choice of the subset of roots. We start with the explicit example of $(1,2)$-flags in $\SO_0(p,3)$, which serves as an introduction to the general methods used in the subsequent Sections. Note that the component count for full flags of $\SO_0(q,q)$ was achieved by Zelevinsky in \cite{Zel00} and the one for full flags in $\SO_0(q+1,q)$ was achieved by Gekhtman--Shapiro--Vainshtein in \cite{GSV03}, both of these result using the formalism of cluster algebras. Our proof is in the same spirit, except that some of the cluster variables would need to be vector-valued for it to make sense.
    
    In the space of full flags of a split-real Lie group (i.e. $\SO_0(q+1,q)$ or $\SO_0(q,q)$), there are special connected components in $\Omega(F_0)\cap\Omega(F_\infty)$ called \emph{totally positive} introduced by Lusztig in \cite{Lus94}. These connected component are made of exactly one cell, hence are homeomorphic to open balls. In the space of $(1,\dots, q-1)$-flags in $\SO_0(p,q)$ with $p\neq q$, there are also special connected components in $\Omega(F_0)\cap\Omega(F_\infty)$ called \emph{$\Theta$-positive} introduced by Guichard--Wienhard in \cite{GW24}. These components are also made of exactly one cell, hence diffeomorphic to open balls.

    \subsection{\texorpdfstring{$(1,2)$}{(1,2)}-flags in \texorpdfstring{$\SO_0(p,3)$}{SO(p,3)}}\label{sec:12flagsSOp3}

    We will now compute the number of connected components of $\Omega(F_0) \cap \Omega(F_{\infty})$ for the space of $(1,2)$-flags in $\SO_0(p,3)$, meaning for the flag varieties of the form $\SO_0(p,3)/P_{1,2}$ for some $p > 3$. \\
    
\begin{figure}[h!]
    \begin{center}
      \includegraphics[width=.25\linewidth]{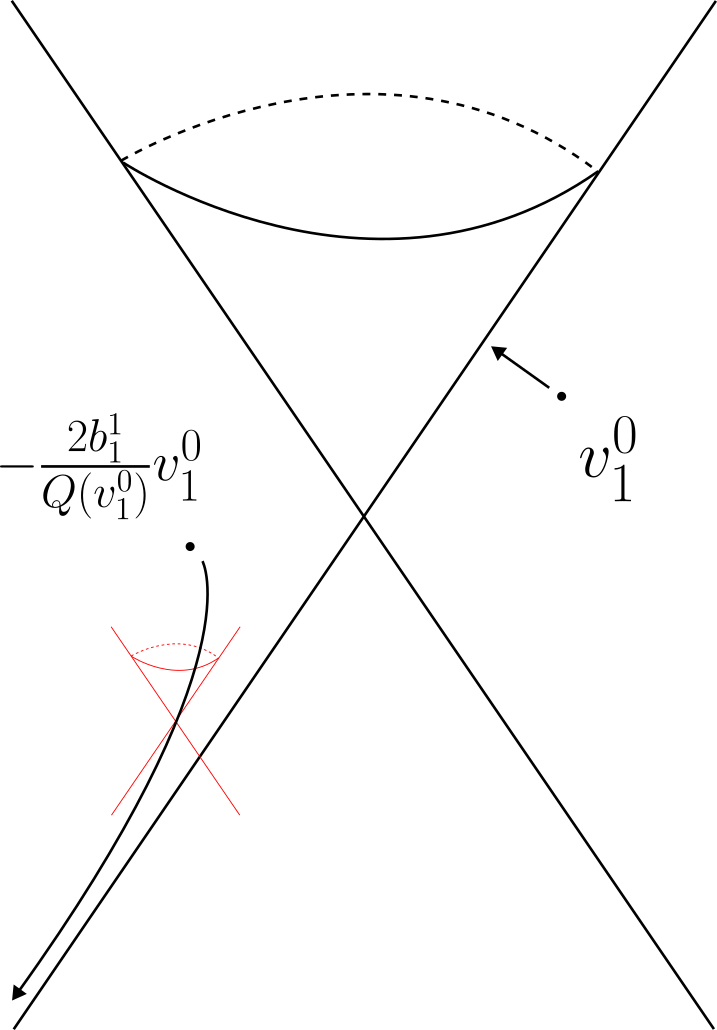}
      \end{center}
      \caption[]{The transition from $v_1^0$ spacelike to $v_1^0$ future type when $b_1^1 > 0$. In red, the light cone of $-\frac{2 b_1^1}{Q(v_1^0)} v_1^0$.}\label{transisiton}
    \end{figure}

    The elements $v_1^0$ and $v_2^0$ both belong to a Lorentzian space $\R^{p-2,1}$. We want now to identify which cell glue to which when changing the sign of $Q(v_1^0)$ in $\Omega(F_0) \cap \Omega(F_{\infty})$. Let's assume for instance that $Q(v_1^0) > 0$ and $b_1^1 > 0$. The point $v_2^0$ must not be on the lightcone of $- \frac{2 b_1^1}{Q(v_1^0)} v_1^0$ (because we have $Q(v_2^{0,(1)})\neq 0$) and is therefore confined into three areas, namely space, future and past. When $v_1^0$ goes from the space part to the future part, $-\frac{2 b_1^1}{Q(v_1^0)} v_1^0$ is negatively colinear to $v_1^0$ and goes to infinity along the past orientation of a photon; therefore the cell corresponding to the past of $-\frac{2 b_1^1}{Q(v_1^0)} v_1^0$ does not glue to any other along the locus $Q(v_1^0)=0$ and $b_1^1 > 0$. This is illustrated in Figure \ref{transisiton}.

\begin{figure}[h!]
    \begin{center}
      \includegraphics[width=.3\linewidth]{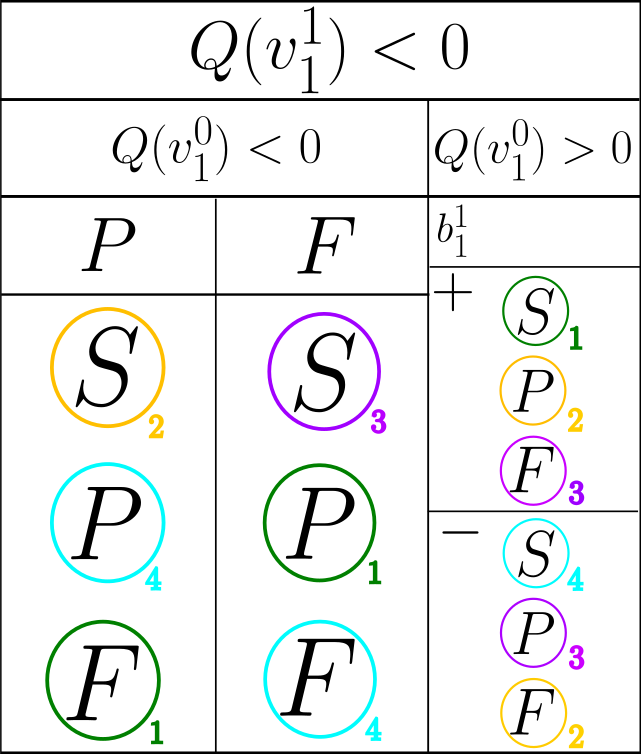}\hspace{2em}\includegraphics[width=.46\linewidth]{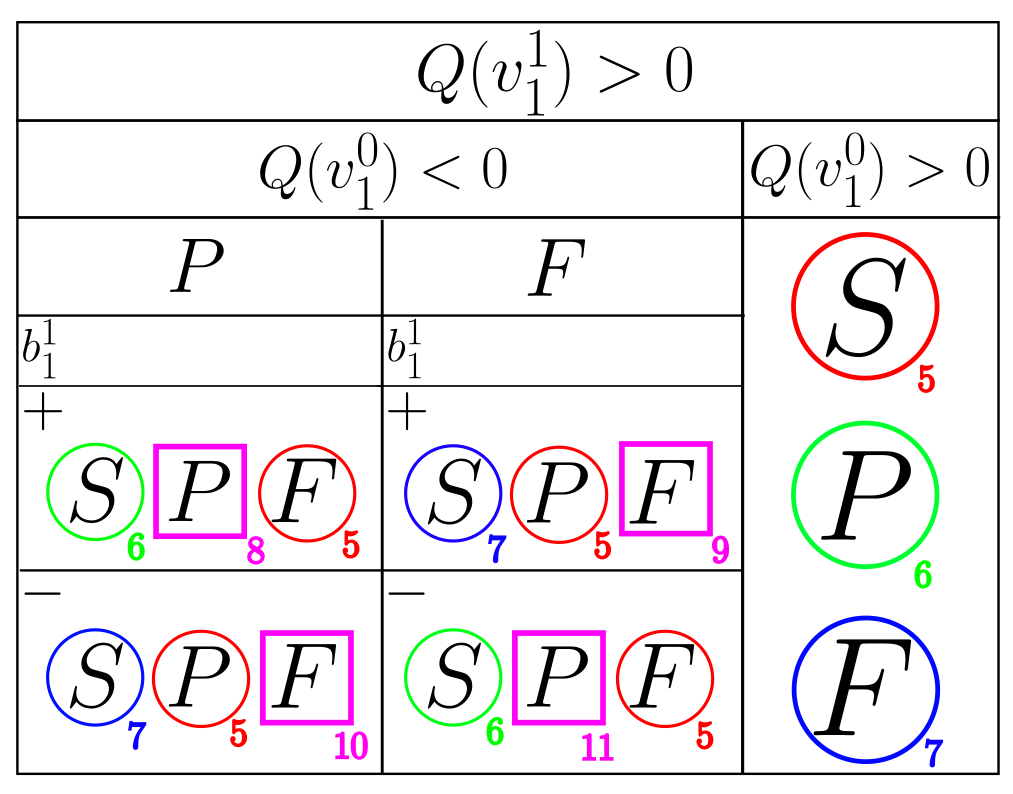}
      \end{center}
      \caption[]{The eleven connected components of $\Omega(F_0) \cap \Omega(F_{\infty})$, one group of four when $Q(v_1^0) < 0$ and one group of seven when $Q(v_1^0) > 0$. The letter S (space), P (past) and F (future) correspond to the cells of $\{Q(v_2^{0, (1)}) \neq 0\}$ when $Q(v_1^0)$ and $Q(v_1^1)$ have constant sign. When $Q(v_1^0)<0$, the column is divided into subcolumns P (where $v_1^0$ is past type) and F (where $v_1^0$ is future type). Cells belonging to the same connected components are circled in the same color and labeled by the same number. The pink components (8,9,10,11) are isolated and are the $\Theta$-positive components.}\label{photonpointe}
    \end{figure}

    Inversely, when going from $v_1^0$ of future type to $v_1^0$ spacelike with $b_1^1$ positive the cell where $v_2^{0,(1)}$ is future type does not glue to any other. When comparing, one then sees that the cell where $v_2^{0,(1)}$ is future type and $v_1^0$ spacelike is glued to the cell where $v_2^{0,(1)}$ is spacelike and $v_1^0$ is future type, and the cell with $v_2^{0,(1)}$ spacelike and $v_1^0$ spacelike is glued with the cell where $v_2^{0,(1)}$ is past like and $v_1^0$ future type. When doing every possible transition, one can then count the exact number of connected components, see Figure \ref{photonpointe}.

    \begin{prop}
        The space $\Omega(F_0) \cap \Omega(F_{\infty})$ in $\SO_0(p,3)/P_{1,2}$ has $11$ connected components, $4$ of which are $\Theta$-positive.
    \end{prop}

    \subsection{\texorpdfstring{$(1,\dots,k)$}{(1,...,k)}-flags in \texorpdfstring{$\SO_0(p,q)$}{SO(p,q)}}

    We will now be dealing with the case of complete flags of isotropic spaces up to dimension $k < q-1$ in $\SO_0(p,q)$ for $p\geq q$. As in the previous sections, we know that when for all $0 \leq m \leq k-1$ and $i < k-j$, $Q(v_j^{i, (m)}) \neq 0$, the equations of transversality become
    
    \[ \begin{split}
        Q(v_1^{k-1, (0)}) &\neq 0, \\
        Q(v_1^{k-2, (0)})Q(v_2^{k-2, (1)}) &\neq 0, \\
        ...& \\
        Q(v_1^{0, (0)})Q(v_2^{0, (1)})...Q(v_k^{0, (k-1)}) &\neq 0. \\
    \end{split}\]
    
    It is easy to check what happens when all $Q(v_j^{i, (j-1)})$ are non-zero except $Q(v_k^{0, (k-1)})$. In this case the equation becomes $Q(v_k^{0, (k-1)}) \neq 0$ and since $v_k^{0, (k-1)}$ is in the pseudo-Riemannian space $\R^{p-k, q-k}$, $v_k^{0, (k-1)}$ can either be spacelike or timelike. This holds for every cell in $\SO_0(p,q)/P_{1,...,k}$, however since the transversality equations do not forbid $Q(v_j^{i, (j-1)}) = 0$ when $i \neq k-j$, we must track which cells glue to which along $Q(v_j^{i, (j-1)}) = 0$ when $i \neq k-j$. We only have to explicit what happens when $Q(v_1^{i, (0)})$ changes sign as when those are constant we can deduce the rules for the $Q(v_j^{i, (j-1)})$, $j > 1$ by induction on $j$. \\

    \begin{lemma}\label{lemma1}
        We have that

        \[a_j^{i, (1)} = a_j^i + \frac{2 b_1^{k-1} a_1^i}{Q(v_1^i)} + \sum_{i+1}^{k-j} \frac{2 a_j^l b_1^l a_1^i}{Q(v_1^i)}, \]
        
        \[v_j^{i, (1)} = v_j^i + \frac{2 b_1^{k-j+1}}{Q(v_1^{i})} v_1^i + \sum_{l=i}^{k-j} \frac{a_j^{l} b_1^l}{Q(v_1^i)} v_1^i.\]
    \end{lemma}

    \begin{proof}
        Let us fix $j$ and prove by decreasing induction on $i$ that

        \[ \frac{2 b_1^{k-1} }{Q(v_1^{k-2})} + \sum_{l=i+1}^{k-j} \frac{2 a_j^{l, (1)} b_1^l }{Q(v_1^{l-1})} = \frac{2 b_1^{k-1} }{Q(v_1^i)} + \sum_{l=i+1}^{k-j} \frac{2 a_j^{l} b_1^l }{Q(v_1^i)}.\]

        The result is true for $i = k-j$ since $a_j^{k-j, (1)} = a_j^{k-j} + \frac{2, b_1^{k-1} a_1^{k-j}}{Q(v_1^{k-j})}$. Assume that it is true for any $i' > i$. Then 

        \[ \begin{split}
            \frac{2 b_1^{k-1} }{Q(v_1^{k-2})} + \sum_{l=i+1}^{k-j} \frac{2 a_j^{l, (1)} b_1^l }{Q(v_1^{l-1})} &= \frac{2 b_1^{k-1} }{Q(v_1^{i+1})} + \sum_{l=i+2}^{k-j} \frac{2 a_j^{l} b_1^l }{Q(v_1^{i+1})} + \frac{2 a_j^{i+1, (1)} b_1^{i+1}}{Q(v_1^i)} \\
            &= \frac{2 a_j^{i+1} b_1^{i+1}}{Q(v_1^i)} + \frac{2 b_1^{k-1} }{Q(v_1^{i+1})} + \frac{4 b_1^{k-1}a_1^{i+1} b_1^{i+1}}{Q(v_1^{i+1})Q(v_1^i)} \\
            &+ \sum_{l=i+2}^{k-j} \frac{2 a_j^{l} b_1^l }{Q(v_1^{i+1})} + \sum_{l=i+2}^{k-j} \frac{4 a_j^lb_1^la_1^{i+1}b_1^{i+1}}{Q(v_1^{i+1})Q(v_1^i)} \\
            &= \frac{2 a_j^{i+1} b_1^{i+1}}{Q(v_1^i)} + \frac{2 b_1^{k-1} }{Q(v_1^{i+1})} \left(\frac{Q(v_1^i) + 2 a_1^{i+1}b_1^{i+1}}{Q(v_1^i)}\right) \\
            &+ \sum_{l=i+2}^{k-j} \frac{2 a_j^{l} b_1^l }{Q(v_1^{i+1})} \left(\frac{Q(v_1^i) + 2 a_1^{i+1}b_1^{i+1}}{Q(v_1^i)} \right) \\
            &= \frac{2 a_j^{i+1} b_1^{i+1}}{Q(v_1^i)} + \frac{2 b_1^{k-1} }{Q(v_1^{i})} + \sum_{l=i+2}^{k-j} \frac{2 a_j^{l} b_1^l }{Q(v_1^{i})} \\
            &= \frac{2 b_1^{k-1} }{Q(v_1^i)} + \sum_{l=i+1}^{k-j} \frac{2 a_j^{l} b_1^l }{Q(v_1^i)},
        \end{split}\]

        using that $Q(v_1^{i+1}) = Q(v_1^i) + 2 a_1^{i+1}b_1^{i+1}$. This gives us the intended results.
    \end{proof}

    \begin{rem}
        In particular, when $Q(v_1^{l,(0)})$ goes to zero for $l \neq i$, $a_j^{i, (1)}$ and $v_j^{i, (1)}$ are bounded.
    \end{rem}

    \begin{lemma}
        When $Q(v_1^{0,(0)})$ goes to zero, we have for all $j \geq 2$,

        \[b_j^{1, (1)} \sim - \frac{Q(v_j^{0, (1)})}{2 a_j^{1,(1)}}.\]
    \end{lemma}

    \begin{proof}
        We know that $Q(v_j^{1, (1)}) = Q(v_j^{0, (1)}) + 2 a_j^{1, (1)} b_j^{1, (1)}$, meaning that 

        \[b_j^{1, (1)} = \frac{Q(v_j^{1, (1)}) - Q(v_j^{0, (1)})}{2 a_j^{1, (1)}}.\]

        From lemma \ref{lemma1} we know that $Q(v_j^{1, (1)})$ and $a_j^{1, (1)}$ are both bounded when $Q(v_1^{0,(0)})$ goes to zero, meaning that we get the equivalence

        \[b_j^{1, (1)} \sim - \frac{Q(v_j^{0, (1)})}{2 a_j^{1,(1)}}.\]
    \end{proof}

    \begin{lemma}\label{Change1}
        For all $i$, when $Q(v_1^{i,(0)})$ goes to zero, $Q(v_j^{l, (j-1)})$ stays bounded for all $(j,l)\neq (2,i)$.
    \end{lemma}

    \begin{proof}
        This lemma relies on similar computations to those presented in the proof of Lemma \ref{lemma1}.
    \end{proof}

    Let us then write down the effect of changing the signs of $Q(v_1^{l,(0)})$.

    \begin{prop}\label{rules}
            When $Q(v_1^{i,(0)})$ goes to zero with $i \geq 1$, 

            \[v_2^{i, (1)} \sim \frac{2 a_2^{i, (1)}b_1^{i}}{Q(v_1^i)} v_1^i,\]

            \[a_2^{i, (1)} \sim \frac{2 a_1^{i+1, (1)}b_1^{i+1, (1)} a_1^i}{Q(v_1^i)},\]

            while all other parameters stay bounded. 
            When $Q(v_1^{0,(0)})$ goes to zero, 

            \[v_2^{0, (1)} \sim \frac{2 a_2^{1, (1)} b_1^{1, (0)}}{Q(v_1^0)} v_1^{0},\]

            \[b_2^{1, (1)} \sim - \frac{ Q(v_2^{0, (1)})}{2 a_2^{1, (1)}}\]

            while all other parameters stay bounded.
    \end{prop}

    This result shows that two different phenomena appear depending of whether $Q(v_1^{0, (0)})$ changes sign or any of the $Q(v_1^{i, (0)})$ for $i \geq 1$ changes sign. Let us explicit what happens in those two cases.

    \begin{lemma}\label{prop:transition_pseudo_riemannian}
        When $Q(v_1^{i, (0)})$ changes sign, each cell where $v_2^{i, (1)}$ is spacelike is glued to the cell where $v_2^{i, (1)}$ is timelike and all other parameters are the same. Furthermore $a_2^{i, (1)}$ must have the sign of $\frac{a_1^{i+1, (1)} b_1^{i+1, (1)} a_1^{i,(0)}}{Q(v_1^{i,(0)})}$ before and after the change of sign of $Q(v_1^{i, (0)})$. 
    \end{lemma}

    \begin{proof}
        Since $Q(v_2^{i, (1)}) \sim \frac{2 a_2^{i, (1)}b_1^{i, (0)}}{Q(v_1^{i, (0)})}$, $v_2^{i, (1)} \neq 0$ implies that $v_2^{1, (0)}$ is outside the isotropic cone of a point $v_2^{1, (0)} - v_2^{i, (1)}$ which goes to infinity along the isotropic line $\Span(v_1^{0, (0)})$ in one direction or the other depending on the sign of $-\frac{a_2^{i, (1)}b_1^{i, (0)}}{Q(v_1^{i, (0)})}$. Since $Q(v_1^{i, (0)})$ changes sign, this implies that the direction of $\Span(v_1^{i, (0)})$ along which $v_2^{i, (1)} - v_2^{i, (0)}$ goes to infinity is reversed when $Q(v_1^{i, (0)})$ goes to zero positively or negatively, which implies that the time and space parts of $Q(v_2^{i, (1)}) \neq 0$ are reversed, hence the result. The second part follows immediately from $a_2^{i, (1)} \sim \frac{2 a_1^{i+1, (1)}b_1^{i+1, (1)} a_1^{i, (0)}}{Q(v_1^{i, (0)})}$. All other parameters remain bounded when $Q(v_1^{i, (0)})$ goes to zero.
    \end{proof}

    \begin{lemma}
        When $Q(v_1^{0, (0)})$ changes sign, each cell where $v_2^{0, (1)}$ is spacelike is glued to the cell where $v_2^{0, (1)}$ is timelike and all other parameters are the same. Furthermore $b_2^{1, (1)}$ must have the sign of $- \frac{ Q(v_2^{0, (1)})}{2 a_2^{1, (1)}}$ before and after $Q(v_1^{0, (0)})$ changes sign.
    \end{lemma}

    \begin{proof}
        Since $k \leq q-2$, the space $\R^{p-k, q-k}$ to which $v_2^{0}$ belongs to is non-Lorentzian pseudo-Riemannian and the same argument as before applies. The second part follows immediately from $b_2^{1, (1)} \sim - \frac{ Q(v_2^{0, (1)})}{2 a_2^{1, (1)}}$. All other parameters remain bounded when $Q(v_1^0)$ goes to zero.
    \end{proof}

    \begin{prop}\label{1kflags}
        The space $\Omega(F_0) \cap \Omega(F_{\infty}) \subset \SO_0(p,q)/P_{1,...,k}$ has $2^k$ connected components, determined by the signs of the upper left minors of size $i\times i$ of the matrix $S$ defined in Section \ref{sec:general_eq_transversality} for $1 \leq i \leq k$.
    \end{prop}

    \begin{proof}
        By applying the previous rules when changing the signs of the $Q(v_j^{i, (j-1)})$ that can vanish with $F$ staying transverse to $F_0$, one can then see that every change only reverses the space and time parts defined by the corresponding equations, and thus that every equation of transversality only split $\SO_0(p,q)/P_{1, ..., k}$ into two, thus giving $2^k$ connected components as there are $k$ equations.
    \end{proof}

    \subsection{\texorpdfstring{$(1,\dots,q-1)$}{(1,...,q-1)}-flags in \texorpdfstring{$\SO_0(p,q)$}{SO(p,q)}}\label{submaximalflags}

    In this section we will only treat the case where $p > q$, the case where $p=q$ will be discussed in a following section. We know from the works of Guichard--Wienhard on $\Theta$-positivity (see \cite{GW24}) that there must be $2^{q-1}$ $\Theta$-positive connected components in $\Omega(F_0) \cap \Omega(F_{\infty})$. We apply the same reasoning as before, except that the elements $v_j^{0, (j-1)}$ now belong to the Lorentzian space $\R^{p-q+1, 1}$, meaning that we now have to distinguish between the space part, future part and past part of the complement of a lightcone. When $Q(v_1^{i, (0)})$ changes sign for $i > 0$, the same happens as in the previous case as $v_2^{i, (1)}$ belongs to a non-Lorentzian pseudo-Riemannian space. We start by describing which cells connect to which along $Q(v_1^{0, (0)})=0$. The following lemma is a reformulation of the computation presented in \ref{sec:12flagsSOp3}:\\

    \begin{lemma}\label{prop:transition_lorentzian}
        When $a_2^{1, (1)} b_1^{1, (0)}$ is positive, the cells with $v_1^{0, (0)}$ spacelike and $v_2^{0, (1)}$ of future type is connected to the ones with $v_1^{0, (0)}$ future type and $v_2^{0, (1)}$ spacelike. The cells with $v_1^{0, (0)}$ spacelike and $v_2^{0, (1)}$ spacelike is connected to the ones with $v_1^{0, (0)}$ future type and $v_2^{0, (1)}$ past type, as well as the cells with $v_1^{0, (0)}$ past type and $v_2^{0, (1)}$ future type. The cells with $v_1^{0, (0)}$ spacelike and $v_2^{0, (1)}$ of past type is connected to the ones with $v_1^{0, (0)}$ past type and $v_2^{0, (1)}$ spacelike. In all the mentioned cells $b_2^{1, (1)}$ must have the sign of $- \frac{ Q(v_2^{0, (1)})}{2 a_2^{1, (1)}}$, and all other parameters remain unchanged.
        
        When $a_2^{1, (1)} b_1^{1, (0)}$ is negative, a statement is obtained by substituting ”past type” for ”future type” (resp. ”future type” for ”past type”) in the above. 
    \end{lemma}


\ytableausetup
{boxsize=1em}

    \begin{rem}
        The transversality equations are $\det_i(S)\neq 0$ for $1\leq i\leq q-1$, so the signs of the $\det_i(S)$ separates $\Omega(F_0) \cap \Omega(F_{\infty})$ into at least $2^{q-1}$ connected components. We will show that for a given choice of the signs of the $\det_i(S)$ there are 2 connected components (except in the case when $\Theta$-positive components arise, which we will also detail).
    \end{rem}

    As previously, we start with an initial set of cells each given by the choice of the sign of all $Q(v_j^{i,(j-1)})$ for $1\leq j\leq q-1$ and $0\leq i\leq k-j$. When $Q(v_j^{i,(j-1)})Q(v_j^{i-1,(j-1)})<0$, we also need to fix the sign of $b_j^{i-1,(j-1)}$. We then glue them together along the locus where exactly one of the $Q(v_j^{i,(j-1)})$ vanish. The locus where two or more vanish at the same time is of codimension at least $2$ hence does not change the number of connected components. Contrary to the previous sections, the combinatorics of the gluing process for $(1,\dots,q-1)$-flags is more involved. To proceed with the study of this gluing process, we first describe a synthetic model describing the situation.

    We represent an cell by a anti-triangular $(q-1)\times (q-1)$ matrix $M = (m_{i,j})_{1\leq i\leq q-1,~0\leq j\leq q-1-i}$ with entries valued in $\mathbb Z/3\mathbb Z$ identified with the set $\lbrace *,+,-\rbrace$ (with $*=0,+=1,-=2$), with every entry having one of the two colors red or blue. We call such a matrix a \emph{sign matrix}. The way to read this data is as follows: the color of $m_{i,j}$ corresponds to the sign of $Q(v_i^{j,(i-1)})$ (blue for negative and red for positive), $m_{i,j} = *$ if and only if the $(i,j)$ and $(i,j-1)$ entries are of the same color (the entry $(i,-1)$ count as red as $Q(v_i^{-1,(i-1)})$ is always positive), and if $m_{i,j}\neq *$ then it corresponds to the sign of $b_i^{j,(i-1)}$ (which has to be fixed in that cell). We now want to describe the effect changing the sign of $Q(v_i^{j,(i-1)})$ has on this matrix. We will call such a modification a \emph{alteration} of $M$, and we call all sign matrices that one can obtain from $M$ by successive alteration the \emph{alteration class} of $M$.

    For each $(i,j)$ such that $i<q-1$ and $j< q-1-i$, we define $$\mu_{i,j} = (-1)^{\varepsilon_1}(-1)^{\varepsilon_2} m_{i,j}m_{i+1,j}m_{i,j+1}m_{i+1,j+1}$$ where $$\varepsilon_1 = \left\lbrace\begin{array}{l}
       1\text{ if }m_{i+1,j}\text{ and } m_{i+1,j+1}\text{ are of the same color}\\
       0 \text{ otherwise}
    \end{array} \right. \text{ and }\varepsilon_2 = \left\lbrace\begin{array}{c}
       1\text{ if }m_{i,j} =*\\
       0 \text{ otherwise}
    \end{array} \right.$$
    When $\mu_{i,j} \neq +$ we say that the $(i,j)$ entry of $M$ is \emph{alterable}. When $m_{i,j}$ is alterable, one can alter it as follows: choose a sign ($+$ or $-$) for every $*$ among the four entries $m_{i,j}, m_{i+1,j}, m_{i,j+1}, m_{i+1,j+1}$ such that after replacing every $*$ by the sign chosen, we get $\mu_{i,j} = -$. Then change the color of $m_{i,j}$ and $m_{i+1,j}$, and replace the entries that were not $*$ at the beginning with $*$. There might be multiple sign choices that result in an admissible alteration, hence a alterable entry may have several possible alterations.
    These rules are obtained by translating the result of Proposition \ref{rules} in this setup, and the correction factor given by $\varepsilon_1$ and $\varepsilon_2$ are due to the presence of $a_i^{j,(i-1)}$ in the formulas whereas we are only working with the $b_i^{j,(i-1)}$ here. 
    To get coherent formulas when $j = q-2-i$ (meaning that the $(i+1,j+1)$ entry necessary to the definitions of $\mu_{i,j}$ and $\varepsilon_1$ does not exist), we need to make the following fix: we set $m_{i+1,q-1-i}$ to be a red $+$ when $m_{i+1,q-2-i}$ is blue and a blue $-$ when $m_{i+1,q-2-i}$ is red.
    Propositions \ref{prop:transition_pseudo_riemannian} and \ref{prop:transition_lorentzian} then reformulate to:
    \begin{proposition}
        Two cells lie in the same connected component of $\Omega(F_0) \cap \Omega(F_{\infty})$ if and only if their corresponding sign matrices are in the same alteraion class.
    \end{proposition}
    The remaining of this section is dedicated to counting the different alteration classes of sign matrices. We start with an explicit example which we will need later as it will arise as a special case to treat separately.

    \begin{ex}\label{exemple_cas_particulier}
    Let us explicit an example. Assume $q=4$. The transversality in $\SO_0(p,4)/P_{1,2,3}$ is given by three equations~:

        \[\begin{split}
            Q(v_1^{2, (0)}) &\neq 0, \\
            Q(v_1^{1, (0)})Q(v_2^{1, (1)}) &\neq 0, \\
            Q(v_1^{0, (0)})Q(v_2^{0, (1)})Q(v_3^{0, (2)}) &\neq 0. \\
        \end{split}\]

    Since $v_1^{2, (0)}$, $v_1^{1, (1)}$ and $v_2^{1, (1)}$ all belong to a non-Lorentzian space, the first two equations split the space $\Omega(F_0) \cap \Omega(F_{\infty})$ in four parts. For each cell, $Q(v_1^{2, (0)})$, $Q(v_1^{1, (0)})$, $Q(v_1^{0, (0)})$, $Q(v_2^{1, (1)})$, $Q(v_2^{0, (1)})$ are all non-zero, the third equation of transversality becomes $Q(v_3^{0, (2)}) \neq 0$. The vector $v_3^{0, (2)}$ belong to a Lorentzian space, meaning it can be either spacelike, future or past type.

    Let us start by considering the matrix $$\begin{ytableau}
 {\color{red}*}& {\color{red}*} &{\color{red}*}\\
\color{blue}+ & \color{red}+\\
\color{blue}+
\end{ytableau}$$ Let us then change the nature of $v_1^{0, (0)}$ from spacelike to future with $b_1^{1, (0)} > 0$ using \ref{prop:transition_lorentzian}. Let us assume that $a_2^{1, (1)}$ is negative. Assuming that $Q(v_2^{0, (1)})$ is negative, $- \frac{Q(v_2^{0, (1)})}{2a_2^{1, (1)}}$ must then be positive, meaning that $b_2^{1, (1)}$ must be negative. Finally, since all other parameters remain unchanged, we get the following alteration $$\begin{ytableau}
 {\color{red}*}& {\color{red}*} &{\color{red}*}\\
\color{blue}+ & \color{red}+\\
\color{blue}+
\end{ytableau}\longrightarrow \begin{ytableau}
 {\color{blue}+}& {\color{red}+} &{\color{red}*}\\
{\color{red}*} & {\color{red}*}\\
\color{blue}+
\end{ytableau}$$ Let us now change $v_1^{0, (0)}$ from future to spacelike with $b_2^{1, (1)} > 0$. By applying the same method, one then get $$\begin{ytableau}
 {\color{blue}+}& {\color{red}+} &{\color{red}*}\\
{\color{red}*} & {\color{red}*}\\
\color{blue}+
\end{ytableau}\longrightarrow \begin{ytableau}
 {\color{red}*}& {\color{red}*} &{\color{red}*}\\
\color{blue}- & \color{red}-\\
\color{blue}+
\end{ytableau}$$

    By the applying the same reasoning, one gets the following alterations $$\begin{ytableau}
 {\color{red}*}& {\color{red}*} &{\color{red}*}\\
\color{blue}- & \color{red}+\\
\color{blue}-
\end{ytableau}\longrightarrow \begin{ytableau}
 {\color{blue}+}& {\color{red}+} &{\color{red}*}\\
{\color{red}*} & {\color{red}*}\\
\color{blue}-
\end{ytableau}\longrightarrow \begin{ytableau}
 {\color{red}*}& {\color{red}*} &{\color{red}*}\\
\color{blue}+ & \color{red}-\\
\color{blue}-
\end{ytableau}$$
We are going to show that these two all these sign matrices are part of the same alteration class. Let us start with the matrix $$\begin{ytableau}
 {\color{red}*}& {\color{red}*} &{\color{red}*}\\
\color{blue}+ & \color{red}-\\
\color{blue}-
\end{ytableau}$$ and change the nature of $v_1^{1, (0)}$ from spacelike to timelike. We know that crossing from $v_1^{1, (0)}$ spacelike to $v_1^{1, (0)}$ timelike changes $v_2^{1, (1)}$ from spacelike to timelike. Let us assume that $b_1^{2, (0)}$ and $b_1^{1, (0)}$ are positive. Since $a_2^{1, (1)}$ must have the sign of $\frac{b_1^{2, (0)} a_1^{1, (0)}}{Q(v_1^{1, (0)})}$ and $b_1^{1, (0)}$ positive implies $a_1^{1, (0)}$ negative, we know that $a_2^{1, (1)}$ has to be positive on departure and negative on arrival. On departure, $Q(v_2^{1, (1)}) > 0$ so when $Q(v_2^{0, (1)})$ is negative, $a_2^{1, (1)}$ positive implies $b_2^{1, (1)}$ negative. This correspond to the following alteration:
$$\begin{ytableau}
 {\color{red}*}& {\color{red}*} &{\color{red}*}\\
\color{blue}+ & \color{red}-\\
\color{blue}-
\end{ytableau}\longrightarrow \begin{ytableau}
 {\color{red}*}& {\color{blue}+} &{\color{red}+}\\
{\color{blue}+} & {\color{blue}*}\\
\color{blue}-
\end{ytableau}$$

We will now change the nature of $v_2^{0, (1)}$ from past to spacelike with $b_2^{1, (1)} > 0$, and then change the nature of $v_2^{0, (1)}$ from spacelike to future. By applying the same rules, we then get the following alterations:
$$\begin{ytableau}
 {\color{red}*}& {\color{blue}+} &{\color{red}+}\\
{\color{blue}+} & {\color{blue}*}\\
\color{blue}-
\end{ytableau}\longrightarrow \begin{ytableau}
 {\color{red}*}& {\color{blue}+} &{\color{red}+}\\
{\color{red}*} & {\color{blue}+}\\
\color{red}*
\end{ytableau}\longrightarrow \begin{ytableau}
 {\color{red}*}& {\color{blue}+} &{\color{red}+}\\
{\color{blue}-} & {\color{blue}*}\\
\color{blue}+
\end{ytableau}$$

Finally, let us change $v_1^{1, (0)}$ from timelike to spacelike. With the same reasoning as before, we get the following alteration~:

$$\begin{ytableau}
 {\color{red}*}& {\color{blue}+} &{\color{red}+}\\
{\color{blue}-} & {\color{blue}*}\\
\color{blue}+
\end{ytableau}\longrightarrow \begin{ytableau}
 {\color{red}*}& {\color{red}*} &{\color{red}*}\\
\color{blue}- & \color{red}-\\
\color{blue}+
\end{ytableau}$$

We have shown that all these cells of $\Omega(F_0) \cap \Omega(F_{\infty})$ were part of the same connected component.
\end{ex}

The cells whose sign matrices do not admit any alteration (hence are alone in their respective alteration classes) correspond to the $\Theta$-positive components introduced by Guichard--Wienhard in \cite{GW24}. Indeed, these cells correspond to the only connected components of $\Omega(F_0) \cap \Omega(F_{\infty})$ whose closure are contained in an affine chart of $\mathcal{F}^{1,\dots,q-1}$ (i.e. which are proper), and we will see in \ref{rem:Theta_positive_cells} that there are $2^{q-1}$ such cells, which is also the number of $\Theta$-positive components. 

\begin{definition}
    If a sign matrix $M$ have no admissible alteration, we call $M$ a $\Theta$-positive sign matrix.
\end{definition}

    \begin{lemma}
        In the alteration class of any matrix $M$ there is a matrix $M'$ such that for all $2\leq i\leq q-1$ and for all $0\leq j\leq q-1-i$, the $(i,j)$ entry of $M'$ is blue when $j$ is even and red when $j$ is odd. We call such a matrix a \emph{striped sign matrix}.
    \end{lemma}

    \begin{rem}
        The color of the first line entries of a striped matrix are determined entirely by the signs of the $\det_i(U)$, hence are invariants of the alteration class of $M$.
    \end{rem}

    \begin{proof}
        We prove the result by induction on the column index. If an entry $(i,0)$ with $i\geq 2$ of the first column is red, it is a $*$ by definition, hence the entry $(i-1,0)$ is alterable. Doing any alteration of $(i-1,0)$ will result in $(i,0)$ being blue. By applying this process starting from bottom to the top, we can change the whole first column to blue except for the entry $(1,0)$. Once the column $j$ is either blue or red except for $(1,j)$, the same argument allow us to color the column $j+1$ with the other color, once again except for $(1,j+1)$.
    \end{proof}

    \begin{rem}\label{rem:twice_alter}
        When $M$ is a striped sign matrix and $(i,j)$ is an alterable entry with $i\geq 2$, the only possible alteration at $(i,j)$ takes all four entries $(i,j)$, $(i+1,j)$, $(i,j+1)$ and $(i+1,j+1)$ to $*$, obtaining a matrix that is no longer striped. Altering again the same entry, one get back a striped matrix, and every possible such alteration result in a striped matrix $M'$ such that among the entries $(i,j)$, $(i+1,j)$, $(i,j+1)$ and $(i+1,j+1)$ of $M'$, an even number of them have their opposite signs compared to $M$. We will use this operation of altering twice to change an even number of signs a lot in the following proofs.
    \end{rem}

    \begin{remark}
        Because of the phenomenon described in Remark \ref{rem:twice_alter}, any two sign matrices $M$ and $M'$ which have the same $*$ entries but whose total number of $-$ entries differ mod 2 can not be in the same alteration class. In particular, there are at least two alteration classes of striped matrix having the same signs for all $Q(v_1^{i})$ for $0\leq i\leq q-2$.
    \end{remark}

    \begin{lemma}\label{lem:normal_form}
        If a striped matrix $M$ is not $\Theta$-positive, there is a alteration-equivalent striped matrix $M' = (m'_{i,j})$ such that:
        $$m'_{i,j} = \left\lbrace \begin{array}{l}
            + \text{ if }j \text{ or } q-1-i \text{ is odd}\\
            - \text{ if }j \text{ and } q-1-i \text{ are even}
        \end{array} \right.$$
        unless for all $0\leq j\leq q-2$, $m_{1,j}=*$, in which case
        $$m'_{i,j} = \left\lbrace \begin{array}{l}
            \pm \text{ if } i=2 \text{ and } j=q-3\\
            + \text{else if }j \text{ or } q-1-i \text{ is odd}\\
            - \text{else if }j \text{ and } q-1-i \text{ are even}
        \end{array} \right.$$
        We call a striped matrix of this form a \emph{normalized striped matrix}.
    \end{lemma}

    \begin{proof}
        We start by showing that we can alter $M$ to a striped matrix whose $(2,0)$ entry is alterable if $q$ is odd and whose $(1,0)$ entry is alterable if $q$ is even. If there are no alterable entry with $i\geq 2$, it means that there is one alterable entry $(1,j)$ on the first line. Altering it twice allow us to change the signs of $(2,j)$ and $(2,j+1)$, thus making $(2,j-1)$ alterable because exactly one of the four signs determining if it is alterable changed, and it was not alterable before by hypothesis. There are two edge cases when we can not apply this result, first is when $j = q-2$ and $m_{1,q-2}=m_{1,q-1} =*$: $$\begin{ytableau}
 \none[...] & {*} & {*}\\
\none[...] & \pm
\end{ytableau}$$
in which case $(1,j-1)$ is also alterable and we can apply the previous construction to $(1,j-1)$. The second edge case is when $q=4$ and the first line is all $*$:
$$\begin{ytableau}
 {\color{red}*}& {\color{red}*} &{\color{red}*}\\
\color{blue}\pm & \color{red}\pm\\
\color{blue}\pm
\end{ytableau}$$
We showed how to get to a normal striped form in Example \ref{exemple_cas_particulier}.

        Now we have a alterable entry $(i,j)$ with $i\geq 2$. We will "move" it to make the upper left corner $(2,0)$ alterable. For this, notice that altering $(i,j)$ twice changing the signs of $(i,j)$ and $(i,j+1)$ make $(i,j-1)$ alterable, and altering $(i,j)$ twice changing the signs of $(i,j)$ and $(i+1,j)$ make $(i-1,j)$ alterable. Iterating this process allow us to alter into a striped matrix whose $(2,0)$ entry is alterable. If $q$ is even and the $(1,0)$ is not alterable yet, we can make it alterable by altering twice $(2,0)$ changing the sign of $(2,0)$ and $(3,0)$.
        
        We want now to have a striped matrix for which every $(i,j)$ entry with $q-i-1$ even and $j$ even is alterable. We start with the one on the top left, which we already made alterable, and we will make the one two cells to the right and the one two cells below also alterable. By iterating this process, we will make all of them alterable. If a cell $(i,j)$ is alterable, we can make $(i+2,j)$ and $(i,j+2)$ (provided they exist) alterable by making $(i+1,j+1)$ alterable: if it is not already the case, we can change signs of $(i,j)$ and $(i+1,j+1)$, thus making $(i+1,j+1)$ alterable. Then, if neither $(i+2,j)$ nor $(i,j+2)$ are alterable, we can make them both alterable by altering $(i+1,j+1)$ twice changing signs of $(i+2,j+1)$ and $(i+1,j+2)$. If both $(i+2,j)$ and $(i,j+2)$ are alterable, there is nothing to do. If only $(i+2,j)$ (resp. $(i,j+2)$) is alterable and if $(i+2,j+2)$ is not out of bounds, we can make $(i,j+2)$ (resp. $(i+2,j)$) alterable by altering twice $(i+1,j+1)$ changing signs of $(i+1,j+2)$ and $(i+2,j+2)$ (resp. $(i+2,j+1)$ and $(i+2,j+2)$). If $(i+2,j+2)$ is out of bounds and only exactly one of $(i+2,j)$ or $(i,j+2)$ is alterable, we just ensure that $(i+2,j)$ is alterable by changing signs of $(i+2,j+1)$ and $(i+1,j+2)$. Applying this process starting with the top-left corner and expanding from it ensures that the only the one in the top right. That entry is either $(2,q-4)$ if $q$ is even or $(1,q-3)$ if $q$ is odd. In the latter case, we actually just need to make the $(2,q-3)$ be a $+$ if the first row is not only made of $*$. For this, we distinguish two cases: either $(1,q-4)$ is a $\pm$ or it is a $*$. When $(1,q-4)$ is a $\pm$, we first make it alterable by changing the signs of $(1,q-5)$ and $(1,q-4)$ (since $(1,q-5)$ is alterable) if it is not. Then changing signs of $(1,q-4)$ and $(2,q-3)$ give us the result. When $(1,q-4)$ is a $*$, notice that $(1,q-3)$ and $(1,q-2)$ must be $\pm$ since $(1,q-3)$ is not alterable. Then by hypothesis there is an entry $(1,j)$ which is not a $*$, and we can assume $j$ to be maximal for this property. Then by altering twice successively $(1,j)$, $(1,j+1),\dots,(1,q-4)$, we can change only the signs of $(1,j)$ and $(2,q-3)$, hence the result. When we need to make $(2,q-3)$ alterable, we use the same arguments to change its sign thus making it alterable.

        At this point, we have a striped matrix whose entries $(i,j)$ with $q-1-i$ and $j$ even are alterable (except for $(2,q-4)$ or $(1,q-3)$ when the first row is all $*$ which we will discuss later). We can alter twice each of these entries to make all of them be $-$ and all of their left, bottom and bottom-left neighbors $+$, thus resulting in a normalized striped matrix. When the first row is all $*$, the only entry that can not be changed this way is the last one of the second row which stays as a sign $\pm$.
    \end{proof}

    Note that a normalized striped matrix has all its entries $(i,j)$ with $i\geq 2$ alterable.

    \begin{lemma}
        In the alteration class of any matrix that is not $\Theta$-positive, there is a unique normalized striped matrix such that the first row is either all red $*$ or a sequence of $\pm$ and $*$ such that only the last sign $\pm$ of the row may be $-$.
    \end{lemma}

    \begin{proof}
        Let $M$ be a normalized striped matrix. Let $j_0$ such that $m_{1,j_0}=-$ and there exists a $j_1>j_0$ such that $m_{1,j_1} = \pm$. We can suppose that $j_1$ is minimal for this property. Then for all $j_0<j<j_1$, we have $m_{1,j} =*$. Then by altering twice $(1,j_0),(1,j_0+1),\dots, (1,j_1-1)$, we can change only the signs of $(1,j_0)$ and $(1,j_1)$. When $j_1 = j_0+1$, it may happen that $(1,j_0)$ is not alterable. If that is the case, altering twice $(2,j_0)$ changing the signs of $(2,j_0)$ and $(2,j_0+1)$ makes $(1,j_0+1)$ alterable, then altering twice $(1,j_0+1)$ changing the signs of $(1,j_0+1)$ and $(1,j_0+2)$ makes $(1,j_0)$ alterable, and finally altering twice again $(2,j_0)$ takes the matrix back to being normalized  striped form, except that now $(1,j_0)$ is alterable and we can proceed.
    \end{proof}

    \begin{remark}
        From the previous lemmas we get that the alteration class of a sign matrix (hence a connected component of $\Omega(F_0) \cap \Omega(F_{\infty})$) that is not $\Theta$-positive is entirely determined by the sequence of blue/red entries of the first line and the parity of the number of $-$ entries of any striped sign matrix in its alteration class.
    \end{remark}

    From these lemmas we obtain the following count of the number of connected components of $\Omega(F_0) \cap \Omega(F_{\infty})$ for all $p$ and $q$:

    \begin{prop}\label{1q-1flags}
    In $\SO(p,q)/P_{1,\dots,q-1}$, for $p>q$:
        \begin{itemize}
            \item When $q = 2$, there are $3$ connected components in $\Omega(F_0) \cap \Omega(F_{\infty})$, $2$ of which are $\Theta$-positive.
            \item When $q=3$, there are $11$ connected components in $\Omega(F_0) \cap \Omega(F_{\infty})$, $4$ of which are positive.
            \item When $q \geq 4$, there are $3 \times 2^{q-1}$ connected components in $\Omega(F_0) \cap \Omega(F_{\infty})$, $2^{q-1}$ of which are $\Theta$-positive.
        \end{itemize}
    \end{prop}


    \begin{remark}\label{rem:Theta_positive_cells}
        The $\Theta$-positive components can be obtained in the following way: the corresponding sign matrix must not contain any $*$, hence is striped, with a first line being also striped, i.e. $Q(v_i^{j,(i-1)})$ is positive when $j$ is odd and negative when $j$ is even for all $1\leq i\leq q-1$. Then, the component is determined by the signs $\pm$ of the first line: indeed, since no entry is alterable, fixing any signs in the first line determines uniquely the $(2,q-3)$ entry by $\mu_{2,q-3} =+$, which in turn determines $(2,q-4)$ etc. By computing line by line starting with the second line, going from right to left on each line, one can find a unique $\Theta$-positive striped matrix having the prescribed first line, thus describing the $2^{q-1}$ $\Theta$-positive components. Computing all the vectors $v_i^{j,(i-1)}$ and checking that no entry of the corresponding sign matrix is alterable allow for an explicit algorithm to decide whether a matrix in $U_\Theta$ is $\Theta$-positive. Also note that with our conventions on the quadratic form $Q$, the sign matrix with only $+$ entries is \emph{not} $\Theta$-positive.
    \end{remark}

    \subsection{\texorpdfstring{$(1,\dots,q)$}{(1,...,q)}-flags in \texorpdfstring{$\SO_0(p,q)$}{SO(p,q)}}

    Assume that $p > q+1$. Let~$F \in \Omega(F_0) \cap \Omega(F_{\infty}) \subset \SO_0(p,q)/P_{1,...,q}$. The first $(q-1)$ equations already separate $\Omega(F_0) \cap \Omega(F_{\infty})$ into $3 \times 2^{q-1}$ groups of connected components. When all relevant elements are non-zero, the equation added by the transversality of the maximal isotropic spaces is $Q(v_q^{0, (q-1)}) \neq 0$. However $v_q^{0, (q-1)}$ belongs in the Euclidean space $\R^{p-q,0}$ which is of dimension at least $2$ since $p > q+1$, meaning that the space $Q(v_q^{0, (q-1)}) \neq 0$ is actually connected. 

    \begin{prop}
        The connected components of $\Omega(F_0) \cap \Omega(F_{\infty}) \subset \SO_0(p,q)/P_{1, ..., q}$ are exactly those coming from the transversality in flags of $\SO_0(p,q)/P_{1,...,q-1}$, meaning that the map induced on connected components by $\pi : \SO_0(p,q)/P_{1,...,q} \rightarrow \SO_0(p,q)/P_{1,...,q-1}$ is bijective.
    \end{prop}

    \begin{proof}
        By combining Proposition \ref{rules} with the fact that the $v_j^{0, (j-1)}$ belong to a connected space, one gets the result.
    \end{proof}

    \subsection{\texorpdfstring{$(1,\dots,q)$}{(1,...,q)}-flags in \texorpdfstring{$\SO_0(q+1,q)$}{SO(q+1,q)}}

    Let us now deal with the case of maximal flags when $p = q+1$. Since the group $\SO_0(q+1,q)$ is a split group for which $P_{1,...,q}$ is a Borel subgroup, we know from \cite{Lus94} that there must be $2^q$ totally positive connected components in $\Omega(F_0) \cap \Omega(F_{\infty})$. This stems from the fact that the $v_j^{0, (j-1)}$ now belong to the Euclidean space $\R$ and that the equation $Q(v_q^{0, (q-1)}) \neq 0$ separates $\R$ into two connected components determined by the sign of $v_q^{0, (q-1)}$.

    \begin{prop}
    In $\SO_0(q+1,q)/P_{1,...,q}$:
        \begin{itemize}
            \item When $q = 2$, there are $8$ connected components in $\Omega(F_0) \cap \Omega(F_{\infty})$, $4$ of which are totally positive.
            \item When $q=3$, there are $30$ connected components in $\Omega(F_0) \cap \Omega(F_{\infty})$, $8$ of which are totally positive.
            \item When $q \geq 4$, there are $(q+5)2^{q-1}$ connected components in $\Omega(F_0) \cap \Omega(F_{\infty})$, $2^q$ of which are totally positive.
        \end{itemize}
    \end{prop}

    \begin{proof}
        Let $(v_1^{q-1, (0)},...,v_q^{0, (0)})$ be an element of $\mathfrak{u}_{\Theta}$. As proven in Section \ref{submaximalflags}, the first $q-1$ equations of transversality already splits $\Omega(F^{1,\dots,q-1}_0) \cap \Omega(F^{1,\dots,q-1}_{\infty})$ into $3 \times 2^{q-1}$ connected components, $2^{q-1}$ of which are $\Theta$-positive. Let us assume that $(v_1^{q-1, (0)},...,v_{q-1}^{1, (0)})$ belong to a component of $\Omega(F^{1,\dots,q-1}_0) \cap \Omega(F^{1,\dots,q-1}_{\infty})$ which is not $\Theta$-positive.

        The last equation of transversality is 

        \[[v_1^{0, (0)}]^2...[v_q^{0, (q-1)}]^2 \neq 0,\]
        with the cells where the change of variables are defined being those where all the $v_j^{0, (j-1)}$ are positive or negative for $j < q$. The sign of $v_q^{0, (q-1)}$ separates each connected component of $\Omega(F^{1,\dots,q-1}_0) \cap \Omega(F^{1,\dots,q-1}_{\infty})$ into two cells in $\Omega(F^{1,\dots,q}_0) \cap \Omega(F^{1,\dots,q}_{\infty})$. Depending on the component in which $(v_1^{0, (0)},...,v_{q-1}^{0, (0)})$ is, the signs of the $b_j^{1, (j-1)}$ may be fixed or not. Assume that it is not fixed and let us change the sign of $v_j^{0, (j-1)}$ from positive to negative. Let us first take $b_j^{0, (j-1)} > 0$. Since $v_{j+1}^{0, (j)} \sim \frac{2 b_j^{1, (j-1)}}{v_j^{0, (j-1)}}$, we have $v_{j+1}^{0, (j)} \rightarrow + \infty$ which means that the space $v_{j+1}^{0, (j)} < 0$ does not cross over. However, when taking $b_j^{0, (j-1)}$ negative, the opposite happens, meaning that in the end both $v_{j+1}^{0, (j)} < 0$ and $v_{j+1}^{0, (j)} > 0$ cross over. When assuming that the sign of $b_j^{0, (j-1)}$ is fixed, one then cannot connect all the parts of the cells into two different connected components.

        However, since it was assumed that $(v_1^{q-1, (0)},...,v_{q-1}^{1, (q-2)})$ was not in a $\Theta$-positive component, it is possible to get into a different part of the same connected component in the $(1,\dots,q-1)$-flags in which the sign of $b_j^{0, (j-1)}$ is no longer fixed. This means that a connected component of $\Omega(F^{1,\dots,q-1}_0) \cap \Omega(F^{1,\dots,q-1}_{\infty})$ splits into two connected components in $\Omega(F^{1,\dots,q}_0) \cap \Omega(F^{1,\dots,q}_{\infty})$ when adding the last equation.

        Let us now assume that $(v_1^{q-1, (0)},...,v_{q-1}^{1, (q-2)})$ is in a $\Theta$-positive component. The sign of the $b_j^{1, (j-1)}$ are now all fixed. One then sees from switching the signs of the $v_j^{0, (j-1)}$ that the last equation split the component of $\Omega(F^{1,\dots,q-1}_0) \cap \Omega(F^{1,\dots,q-1}_{\infty})$ into $q+1$ connected components in $\Omega(F^{1,\dots,q}_0) \cap \Omega(F^{1,\dots,q}_{\infty})$, two of which are totally positive (see Figure \ref{fig:maximaux} for an explicit computation in $\SO_0(6,5)$). Combining those two cases gives the result.

        \begin{figure}[hbt]
    \begin{center}
      \includegraphics[width=1.0\linewidth]{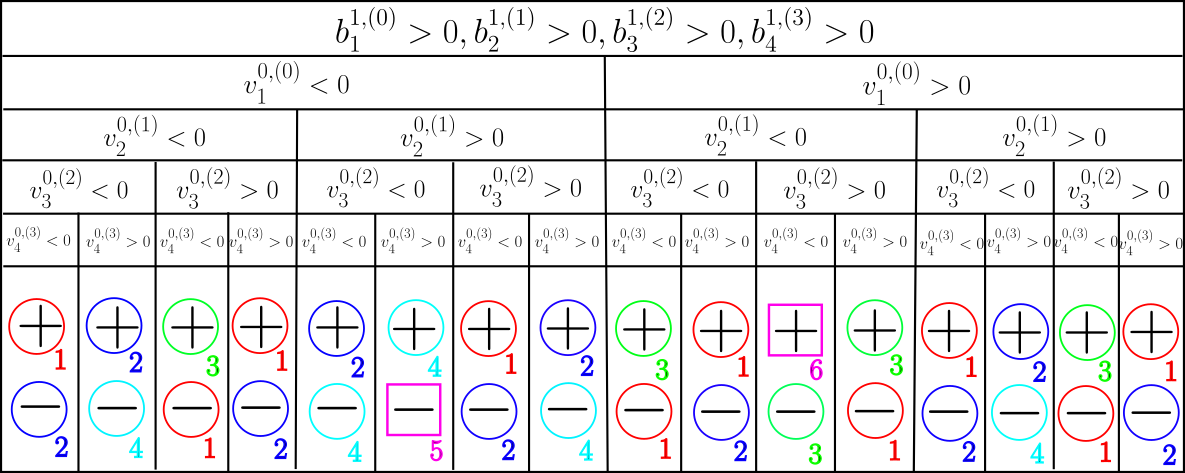}\label{fig:maximaux}
      \end{center}
      \caption[]{The six connected component in a $\Theta$-positive component of the $(1,\dots,4)$-flags in $\SO_0(6,5)$. The signs $+$ and $-$ correspond to the sign of $\{(v_5^{0, (4)})^2 \neq 0\}$ defining a cell where all other relevant values have sign indicated by the columns headers. The components in pink (5,6) are totally positive.}
    \end{figure}
    \end{proof}

    \begin{rem}
        This result was already known from work Gekhtman--Shapiro--Vainshtein, see \cite{GSV03}.
    \end{rem}

    \begin{rem}\label{functionf}
        When all necessary vectors are non-zero, the last equation of transversality becomes

        \[[v_1^{0, (0)}]^2 ... [v_q^{0, (q-1)}]^2 \neq 0.\]

        Note that $[v_1^{0, (0)}]^2 ... [v_q^{0, (q-1)}]^2 = [v_1^{0, (0)} ... v_q^{0, (q-1)}]^2$. Let us define the map $f : \mathfrak{u}_{\Theta} \rightarrow \R$ as

        \[f(v_1^{q-1},...,v_q^{0}) \longmapsto v_1^{0, (0)} ... v_q^{0, (q-1)}.\]

        One may check that $f$ is actually polynomial and thus may be defined on the whole $\mathfrak{u}_{\Theta}$, regardless of the fact that the change of variables may sometimes be undefined. One then gets that $\det(S) = f^2$, and that the sign of $f$ splits $\Omega(F_0) \cap \Omega(F_{\infty})$ into two groups of connected components. In fact, this precisely gives the two connected components composing the groups defined by the first $(q-1)$ equations of transversality which are not $\Theta$-positive.
    \end{rem}
    
    \begin{rem}
        In $\SO_0(3,2)$, this count can be seen geometrically as $\Ein_{2,1}$ is a three dimensional space. Let us write $F_0 = (x_0, \Delta_0)$ and $F_{\infty} = (x_{\infty}, \Delta_{\infty})$ and consider the affine chart from $x_{\infty}$ where $x_0$ is the origin. The space of pointed photons for which the point is in the affine chart is a trivial bundle $\R^{2,1} \times \SS^1$ over $\R^{2,1}$ for which the fiber over $x \in \R^{2,1}$ is the set of photons going through $x$ which is $\Ein_{1,0} \simeq \SS^1$. One may show that in this affine chart, the set of pointed photons transverse to $(x_{\infty}, \Delta_{\infty})$ is $\R^{2,1} \times (\SS^1 \setminus \{N\})$, where $N$ is a point of $\SS^1$.

\begin{figure}[h!]
    \begin{center}
      \includegraphics[width=.5\linewidth]{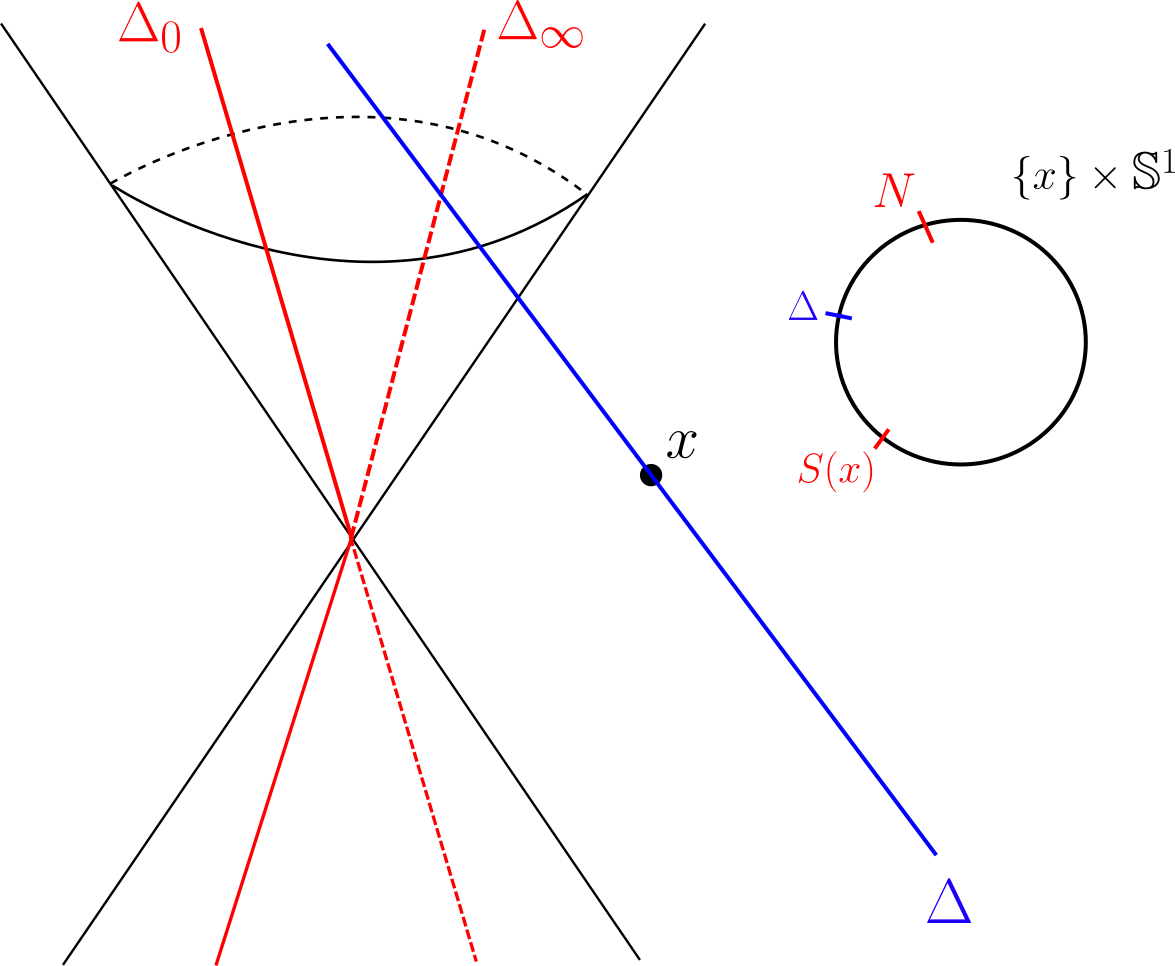}
      \end{center}
      \caption[]{The fiber bundle of pointed photons when $Q(x) > 0$.}
    \end{figure}

    For any point $x \in \R^{2,1} \setminus \Delta_0$, there exists exactly one photon through $x$ intersecting $\Delta_{\infty}$ which is $(x, N)$ and one photon through $x$ intersecting $\Delta_0$ which is $(x, S(x))$. The coordinates of this second photon in $\SS^1$ depend on $x$; those correspond to $v_0^{0, (1)}$ in the affine chart of $\SS^1$ defined by $N$. One must then follow the variations of $S(x)$ when $x$ varies in $\R^{2,1}$. Let $\Delta$ be the only photon through the origin intersecting $\Delta_{\infty}$, i.e the photon of the pointed photon $(0, N)$, and let $V = \Delta_0 \oplus \Delta$. One may show that the two points $S(x)$ and $N$ in $\SS^1$ are equal if and only if $x \in V$. This makes sense as in this case, $v_1^{0, (0)}$ goes to zero and the variable $v_0^{0, (1)}$ is no longer defined.
    
    We can now count geometrically the number of connected components of $\Omega(F_{\infty}) \cap \Omega(F_0)$ in $\SO_0(3,2)/P_{1,2}$. The transversality of $x$ with $x_0$ already separates the connected components into three groups.

    Let's assume that $Q(x) > 0$. When $x$ is in the orthogonal of $V$, $\Delta_0$ and $x + \Delta_0$ have the same endpoint at infinity, thus $S(x) = -N$. When $x \notin V$, $S(x) \neq N$ and the condition of transversality to both $\Delta_0$ and $\Delta_{\infty}$ separates the fiber $\{x\} \times \SS^1$ in two connected components. The vector space $V$ is of dimension $2$ and thus separates $\R^{2,1}$ in two parts. Let $x = (a,b,v)$ be the decomposition of $x$ in $\R^{2,1} = \Delta_{\infty} \oplus \Delta_0 \oplus V^{\perp}$. When $Q(x) > 0$, $x$ may cross the vector space $V$ with both $b > 0$ or $b<0$. Doing so one way or the other makes $S(x)$ go to $N$ clockwise or counterclockwise. Adding all of this up, one eventually counts two connected components, one which contains all the pointed photons with point in $V$ and $b > 0$ and the other which contains those with point in $V$ and $b<0$.

    Assume now that $x$ lies in the future cone of the origin. The same reasoning applies, but we must then always have $b > 0$. When $x$ is on one side of $V$, the fiber $\{x\} \times \SS^1$ separates into two connected components, one of which disappears when $x$ crosses $V$. We must then have three connected components, two of whom correspond to the totally positive components. The same reasoning applies when $x$ lies in the past cone, hence the result.
    \end{rem}

    \subsection{Remaining flags of \texorpdfstring{$\SO_0(p,q)$}{SO(p,q)}}

    Assume that $p \geq q+1$. From the computations we have done we can deduce the number of connected components of $\Omega(F_0) \cap \Omega(F_{\infty}) \subset \SO_0(p,q)/P_{\Theta}$ where $\Theta$ is a subset of the set of roots of $\SO_0(p,q)$ which has not been already covered and $p > q$. In order to do this, one only has to take the description of connected components in $\SO_0(p,q)/P_{1,...,k}$ where $k = \max(\Theta)$ and identify the connected components which are no longer separated by an equation $Q(v_j^{k-j, (j-1)}) \neq 0$ for $j\notin \Theta$. 

    \begin{prop}
    In $\Omega(F_0) \cap \Omega(F_{\infty})$:
        \begin{itemize}
            \item When $\Theta$ does not contain $q$ and is not equal to $\{1,...,q-1\}$, there are $2^{|\Theta|}$ connected components determined by the signs of the significant minors of $S$.
            \item When $\Theta$ contains $q$, is not equal to $\{1,...,q\}$ and $p > q+1$, there are $2^{|\Theta|-1}$ connected components determined by the signs of the significant minors of $S$ except the determinant of $S$ itself which is always non-negative.
            \item When $\Theta$ contains $q$, is not equal to $\{1,...,q\}$ and $p = q+1$, there are $2^{|\Theta|}$ connected components determined by the signs of the significant minor of $S$ except the determinant of $S$ itself which is always non-negative but split in two connected components nonetheless (see Remark \ref{functionf}).
        \end{itemize}
    \end{prop}

    \begin{proof}
        Let us first suppose that $\Theta$ does not contain the last root and $p>q$. Then $\Theta\subset \{1,\dots,q-1\}$ and the projection from $\Omega(F_0^{1,\dots,q-1})\cap \Omega(F_\infty^{1,\dots,q-1})$ to $\Omega(F_0^{\Theta})\cap\Omega(F_\infty^{\Theta})$ obtained by forgetting the subspaces whose dimensions do not belong to $\Theta$ is surjective. However since $\Theta\neq \{1,\dots,q-1\}$ there is $j_0\in \{1,\dots,q-1\}\backslash\Theta$, and the alteration classes of matrices defined in \ref{lem:normal_form} are identified if they differ only from their $j_0$ column. This identifies the two different classes of normalized striped matrices having the same signs of all $\det_i(U)$, showing that there are at most $2^{|\Theta|}$ connected components in $\Omega(F_0^{\Theta}) \cap \Omega(F_{\infty}^{\Theta})$. However since there are still $|\Theta|$ minors that have to be non-zero, there are at least $2^{|\Theta|}$ connected components, hence the result.
    \end{proof}

    \subsection{\texorpdfstring{$q$}{q}-flags in \texorpdfstring{$\SO_0(q,q)$}{SO(q,q)}}\label{maxSOqq}

    For reasons discussed in Section \ref{sec:flag_var_qq}, the case $p=q$ requires a more careful approach. The three main cases which differ from the case $p>q$ are those of $q$-flags, $q'$-flags and $q-1$-flags, which we will discuss in priority.

    \begin{prop}\label{prop:countSOqq}
        When $q$ is even, the subset $\Omega(F_0) \cap \Omega(F_{\infty})$ in $\SO_0(q,q)/P_q$ or $\SO_0(q,q)/P_{q'}$ has two connected components.
    \end{prop}

    \begin{proof}
        As before, the equation of transversality is given by the determinant of the matrix $S$, $\det(S) \neq 0$. However in $\SO_0(q,q)/P_q$, all the $a_j^{i}$ and $v_j^{i}$ are equal to zero and the matrix $S$ is skew-symmetric. In particular, when $q$ is even, $\det(S)$ is always non-negative. It is then well known that the space $\det(S) \neq 0$ in the space of skew-symmetric matrices of even size is the union of two connected components determined by the sign of the Pfaffian of $S$, hence the result. Since by Proposition \ref{transpreserving} there is a transversality-preserving diffeomorphism between the two spaces $\SO_0(q,q)/P_q$ and $\SO_0(q,q)/P_{q'}$, the result also holds in $\SO_0(q,q)/P_{q'}$
    \end{proof}

    \subsection{\texorpdfstring{$(q-1)$}{(q-1)}-flags in \texorpdfstring{$\SO_0(q,q)$}{SO(q,q)}}\label{submaxSOqq}
    
\begin{figure}[h!]
    \begin{center}
      \includegraphics[width=.3\linewidth]{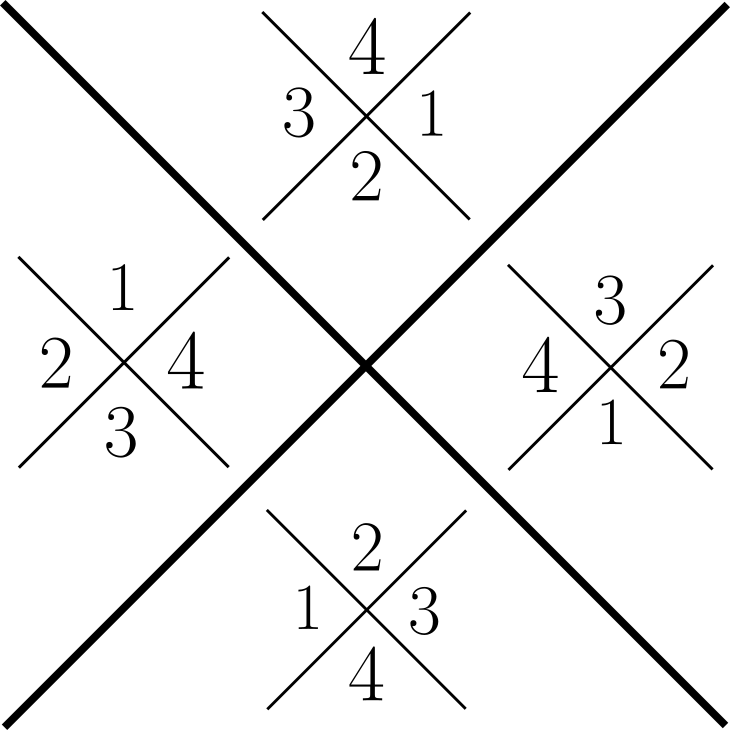}\label{transitionSOqq}
      \end{center}
      \caption[]{In this picture, the big cross represents the space $\{Q \neq 0\}$ in the copy of $\R^{1,1}$ containing $v_1^{0, (0)}$. Each small cross is the space $\{Q \neq 0\}$ in the copy of $\R^{1,1}$ containing $v_2^{0, (1)}$, with the numbers from $1$ to $4$ representing four connected components obtained after identifying which cells glue together.}
    \end{figure}

    The space of $(q-2)$-photons (i.e. $(q-1)$-flags) in $\SO_0(q,q)$ is the flag variety $\SO_0(q,q)/P_{q,q'}$ associated to the last two roots of $\SO_0(q,q)$. It is always self-opposite and the only equation of transversality is given by
    
    \[Q(v_1^{0, (0)})...Q(v_{q-1}^{0, (q-2)}) \neq 0,\]
    
    \noindent where all $v_i^{0, (i-1)}$ belong to $\R^{1,1}$. The space $\{Q \neq 0\}$ in $\R^{1,1}$ is the union of four connected components, which differs from the case $\R^{p,1}$ with $p > 1$ in which the space of spacelike vectors is connected. When $Q(v_i^{0, (i-1)}) \neq 0$ for $i < q-1$, the equation of transversality splits the copy of $\R^{1,1}$ containing $v_{q-1}^{0, (q-2)}$ in four cells given by $Q \neq 0$. Lemma \ref{Change1} and Proposition \ref{rules} tell us that the only transition between cells for which $v_{i+1}^{0, (i)}$ is in different connected component in $\R^{1,1}$ is when $Q(v_{i}^{0, (i-1)})$ changes sign. In this case, Proposition \ref{rules} tells us what happens when $v_1^{0, (0)}$ goes from one connected component of $\{Q(v_1^{0, (0)}) \neq 0\}$ to another. This depends on the sign of $a_2^{0, (1)} b_1^0$; when the value is positive, only two of the four connected components of $\{Q(v_2^{0, (1)}) \neq 0\}$ cross over, while the other two only cross-over while the value of $a_2^{0, (1)} b_{1}^{0}$ is negative. However, since there are no equations regarding the $v_i^{j, (m)}$ for $j > 0$, the signs of $a_2^{0, (1)} b_1^0$ are never fixed in a cell and can be taken to be positive or negative. This allows us to compute what happens to the cell for each transition. The same is true for $i > 1$ by induction. The explicit transitions are illustrated in Figure \ref{transitionSOqq}.
        
    From this, we deduce the following count~:
    
    \begin{prop}
        The subset $\Omega(F_0) \cap \Omega(F_{\infty})$ of $\SO_0(q,q)/P_{q,q'}$ has four connected components.
    \end{prop}
    
    \begin{proof}
        By doing the previous changes, one can take any element of $\Omega(F_0) \cap \Omega(F_{\infty})$ to an element verifying that for all $i < q-1$, $v_i^{0, (i-1)}$ is in the upper connected component of $\{Q \neq 0\}$ in $\R^{1,1}$. The connected component of $\Omega(F_0) \cap \Omega(F_{\infty})$ is then determined by the connected component of $\{Q \neq 0\}$ in which $v_{q-1}^{0, (q-2)}$ belongs, thus giving us four connected components.
    \end{proof}

    \subsection{\texorpdfstring{$(1,\dots,q-1)$}{(1,...,q-1)}-flags in \texorpdfstring{$\SO_0(q,q)$}{SO(q,q)}}
    
    The space of $(1, \dots, q-1)$-flags in $\SO_0(q,q)$ corresponds to the space of full flags of $\SO_0(q,q)$, with $P_{1, \dots, q-1} = P_{1, \dots, q, q'}$ as a Borel subgroup. As such, this count was already achieved by Zelevinsky in \cite{Zel00}. We will only explains how to get this result using our methods without giving a proof.
    
    The situation is completely analogous to Section \ref{submaximalflags}, except for the fact that the $v_i^{0, (i-1)}$ belong to $\R^{1,1}$ in which the space $\{Q \neq 0\}$ has four connected components. The first $q-2$ equations of transversality split $\Omega(F_0) \cap \Omega(F_{\infty})$ in $2^{q-2}$ groups of connected components. In each of those groups, a combinatorial phenomenon similar to the one in Section \ref{submaximalflags} happens, which yields $8$ connected components in each groups when $q > 3$. The cases $q=2$ and $q=3$ are special cases in this regard. Finally, $2^q$ exceptional connected components which are totally positive appear. This gives us the following count~:
    
    \begin{prop}[Zelevinsky, \cite{Zel00}]
    In $\SO_0(q,q)/P_{1,\dots, q-1}$, the number of connected components of $\Omega(F_0) \cap \Omega(F_{\infty})$ is as follows~:
        \begin{itemize}
            \item When $q=2$, there are $4$ connected components, all of which are totally positive.
            \item When $q=3$, there are $20$ connected components, $8$ of which are totally positive.
            \item When $q > 3$, there are $3 \times 2^q$ connected components, $2^q$ of which are totally positive.
        \end{itemize}
    \end{prop}

    \subsection{Remaining flags in \texorpdfstring{$\SO_0(q,q)$}{SO(q,q)}}

    For the same reasons as for $p > q$, the computation for the remaining flags varieties is simply a combination of the previous sections.
    
    \begin{prop}
        Let $\Theta$ be a proper subset of $\{1, \dots q, q'\}$. If $|\Theta \cap \{q, q'\}| = 1$, assume $q$ is even. Then the set $\Omega(F_0) \cap \Omega(F_{\infty})$ in $\SO_0(q,q)/P_{\Theta}$ has $2^{|\Theta|}$ connected components.
    \end{prop}
    
    \begin{proof}
        Any of the first $q-2$ roots in $\Theta$ multiplies the number of connected components by $2$. If $|\Theta \cap \{q,q'\}| = 1$ and $q$ is even, this number is multiplied by $2$ once more by Section \ref{maxSOqq}, thus giving $2^{|\Theta|}$ connected components. If $q,q'$ are both in $\Theta$, this number is multiplied by $4$ by Section \ref{submaxSOqq}, thus also giving $2^{|\Theta|}$ connected components. This yields the result.
    \end{proof}

    \section{Obstructions on the structure of \texorpdfstring{$P_\Theta$}{Theta}-Anosov subgroups of \texorpdfstring{$\SO_0(p,q)$}{SO(p,q)}}

    The aim of this section is to answer a variant of a question by Sambarino regarding the structure of Anosov subgroups.

    \begin{question}\label{question}
        Let $\Gamma$ be a Borel-Anosov subgroup of $\SL(d,\R)$. Is $\Gamma$ virtually isomorphic to either a free group of a surface group ?
    \end{question}

    In \cite{dey2023restrictionsanosovsubgroupssp2nr}, Dey--Greenberg--Riestenberg generalize the arguments developed in \cite{dey2024borelanosovsubgroupsrm} to any $P_{\Theta}$-Anosov subgroup of a semi-simple Lie group $G$ as long as $P_{\Theta}$ is self-opposite. Let $P_{\Theta}$ be a parabolic subgroup of a semi-simple group $G$ such that $G/P_{\Theta}$ is self-opposite. Let $F_0$, $F_{\infty}$ be two transverse points in $G/P_{\Theta}$. The unipotent radical $U_{\Theta}$ acts freely and transitively on $\Omega(F_{\infty})$, which gives a parametrization $g \in U_{\Theta} \mapsto g \cdot F_0$. This parametrization of $\Omega(F_{\infty})$ endows it with a continuous involution $i$ inherited from $U_{\Theta}$ defined by $i(g \cdot F_0) = g^{-1} \cdot F_0$. Since $F_0$ is fixed by $i$ and $i$ preserves transversality with $F_0$, the involution acts on $\Omega(F_0) \cap \Omega(F_{\infty})$ and on $\pi_0(\Omega(F_0) \cap \Omega(F_{\infty}))$.

    \begin{teo}[Dey--Greenberg--Riestenberg, \cite{dey2023restrictionsanosovsubgroupssp2nr}]\label{DGR2}
        Assume that $$i : \pi_0(\Omega(F_0) \cap \Omega(F_{\infty})) \rightarrow \pi_0(\Omega(F_0) \cap \Omega(F_{\infty}))$$ does not have any fixed point. Then any $P_{\Theta}$-Anosov subgroup of $G$ is virtually isomorphic to either a free group or a surface group.
    \end{teo}


    Our aim is to compute the action of the involution on the connected components of $\Omega(F_0) \cap \Omega(F_{\infty})$ for any parabolic subgroups $P_{\Theta}$ of $\SO_0(p,q)$, to answer positively to Question \ref{question} (when replacing $\SL(d,\mathbb{R})$ with $\SO(p,q)$) when possible.

    \subsection{Computing the involution}

    In order to parametrize the space $\Omega(F_{\infty})$ in $\SO_0(p,q)/P_{1,...,k}$, we used the diffeomorphism from $\mathfrak{u}_{\Theta}$ to $U_{\Theta}$,

    \[(v_1^{k-1},...,v_k^{0}) \in \mathfrak{u}_{\Theta} \longmapsto \exp(u_1(v_1^{k-1})) ... \exp(u_k(v_k^{0})) \cdot F_0 \in \SO_0(p,q)/P_{1,...,k}.\]

    We wish to compute the effect of the involution on the coordinates in $\mathfrak{u}_{\Theta}$. For each $(v_1^{k-1},...,v_k^{0})$ in $\mathfrak{u}_{\Theta}$, there exists an element $(i(v_1^{k-1}),...,i(v_k^{0}))$ in $\mathfrak{u}_{\Theta}$ such that 

    \[\begin{split}
        [\exp(u_1(v_1^{k-1})) ... \exp(u_k(v_k^{0}))]^{-1} &= \exp(- u_k(v_k^{0}))...\exp(- u_1(v_1^{k-1})) \\
        &= \exp(u_1(i(v_1^{k-1}))) ... \exp(u_k(i(v_k^{0}))).
    \end{split}\]

    We will write 

    \[\begin{split}
        i(v_1^{k-1}) &= (i(a_1^{k-1}),...,i(a_1^{1}),i(v_1^0),i(b_1^1),...,i(b_1^{k-1})), \\
        &... \\
        i(v_{k-1}^{1}) &= (i(a_{k-1}^{1}),i(v_{k-1}^0),i(b_{k-1}^1)), \\
        i(v_k^0) &= i(v_k^0).
    \end{split}\]

    \begin{prop}
        The involution on $\mathfrak{u}_{\Theta}$ is as follows~:

        \[ \begin{split}
            i(a_j^i) &= - a_j^i + \sum_{l=1}^{k-i-j} i(a_j^{i+l}) i(a_{j+l}^i), \\
            i(v_j^0) &= -v_j^0 + \sum_{l=1}^{k-j} i(a_{j+l}^{1}) i(v_{j+l}^0), \\
            i(b_j^{i}) &= - b_j^{i} - Q(v_i^{i}, v_j^{i}) - \sum_{l=1}^{k-j-i+1} a_j^{k-i-l+1} i(b_{j+l}^i).
        \end{split}\]
    \end{prop}

    \begin{proof}
        This is a straightforward computation. Of particular note is that the involution acts by transposition on the matrix $S$, providing further evidence that $i$ preserves the transversality to $F_0$.
    \end{proof}

    \subsection{\texorpdfstring{$(1,\dots,k)$}{(1,...,k)}-flags in \texorpdfstring{$\SO_0(p,q)$}{SO(p,q)}}

    As explained in Proposition \ref{1kflags}, the connected components of $\Omega(F_0) \cap \Omega(F_{\infty})$ in $\SO_0(p,q)/P_{1,...,k}, k <q-1$ are determined by the signs of the successive upper-left minors of $S$. Since the involution acts on $S$ by transposition, those minors are left unchanged and every connected component of $\Omega(F_0) \cap \Omega(F_{\infty})$ is preserved by the involution.

    \begin{prop}
        When $\Theta = \{1,...,k\}$ and $k <q-1$, the action of $i$ on $\pi_0(\Omega(F_0) \cap \Omega(F_{\infty}))$ is trivial.
    \end{prop}

    \subsection{\texorpdfstring{$(1,\dots,q-1)$}{(1,...,q-1)}-flags in \texorpdfstring{$\SO_0(p,q)$}{SO(p,q)}}

    In order to compute the action of $i$ on the connected components of $\Omega(F_0) \cap \Omega(F_{\infty})$, one only has to exhibit one element of each connected component on which the involution is easily computable. Let us first assume that for each $i,j$, $a_j^{i} = 0$. We deduce from the formulas of the change of variables that for each $i,j,m$, $a_j^{i, (m)} = 0$. This implies that for each $j$, $Q(v_j^{k-j, (j-1)}) = Q(v_j^{k-j-1, (j-1)}) = ... = Q(v_j^{0, (j-1)})$. Since only the signs of $Q(v_1^{k-1, (0)}), Q(v_2^{k-2, (1)}),...,Q(v_{q-2}^{0, (q-3)})$ are fixed, it is possible to find a cell of any connected components determined by the first $q-2$ equations in which all the $a_j^{i}$ are equal to zero. Since the $\Theta$-positive connected components are never stable by the involution, we will not discuss them.

    As the involution preserves the signs of the upper-left minors, and since all the $a_j$ are equal to zero, it must also preserve each equation set ${Q(v_j^{0, (j-1)}) \neq 0}$. Let us determine the action of the involution on its connected components.
    
    \begin{prop}\label{involutioncompute}
        For each $j$, the action of $i$ on ${Q(v_j^{0, (j-1)}) \neq 0}$ always sends the space component on itself. When $j$ is even, the future and past components are preserved. When $j$ is odd, the future and past components are exchanged.
    \end{prop}

    \begin{proof}
        For $j=1$, the result is clear since $i(v_1^0) = - v_1^0$. For $j=2$, we have $v_2^{0, (1)} = v_2^0 + \frac{2 b_1^{q-2}}{Q(v_1^{0})} v_1^0$. Let us take $v_2^0 = 0$. Since the involution acts by $i(b_1^{q-2}) = - b_1^{q-2} - Q(v_1^0, v_2^0)$, we now have $i(b_1^{q-2}) = - b_1^{q-2}$. Since $i(v_1^0) = -v_1^0$, this gives us $i(v_2^{0, (1)}) = v_2^{0, (1)}$, hence the result.

        Let us take $v_{\ell}^0 = 0$ and $b_{\ell}^{i} = 0$ for each $\ell$ even. For $j=3$, we have 

        \[v_3^{0, (2)} = v_3^0 + \frac{2 b_1^{q-3}}{Q(v_1^0)}v_1^0 + \frac{2 b_2^{q-3, (1)}}{Q(v_2^{0, (1)})}v_2^{0, (1)}.\]

        Since $v_2^0 = 0$ and $b_2^{q-3} = 0$, we have 

        \[b_2^{q-3, (1)} = - \frac{2 b_1^{q-2}b_1^{q-3}}{Q(v_1^0)},\]
        as well as $v_2^{0, (1)} = \frac{2 b_1^{q-2}}{Q(v_1^0)}v_1^0$. This gives us

        \[\begin{split}
            v_3^{0, (2)} &= v_3^0 + \frac{2 b_1^{q-3}}{Q(v_1^0)} - \frac{2 b_1^{q-3}}{Q(v_1^0)} \\
            &= v_3^0.
        \end{split}\]

        By repeating this process, we obtain that when $j$ is even, $v_j^{0, (j-1)} = \frac{2 b_{j-1}^{k-j+1}}{Q(v_{j-1}^{0, (j-2)})} v_{j-1}^{0, (j-2)}$ and when $j$ is odd, $v_j^{0, (j-1)} = v_j^0$. This gives us the intended result.
    \end{proof}

    The first $q-2$ equations of transversality split $\Omega(F_0) \cap \Omega(F_{\infty})$ into $2^{q-2}$ groups of connected components. When all relevant variables are non-zero, the last equation of transversality is

    \[Q(v_1^{0, (0)})Q(v_2^{0, (1)}),...,Q(v_{q-1}^{0, (q-2)}) \neq 0.\]

    The sign of the determinant of $S$ once again splits each of the $2^{q-2}$ groups into two groups. Finally, as proven in \ref{1q-1flags}, the parity of the number of timelike vectors in the $v_1^{0, (0)},...,v_{q-1}^{0, (q-2)}$ which are in the future cone splits each of these groups into two connected components.

    \begin{prop}
        Among the $3 \times 2^{q-1}$ connected components in $\Omega(F_0) \cap \Omega(F_{\infty})$, the $2^{q-1}$ $\Theta$-positive components are all exchanged, $2^{q-1}$ of the remaining components are exchanged and the rest are stable by $i$.
    \end{prop}

    \begin{proof}
        Assume that $q$ is even. The involution thus acts on $\{Q(v_{q-1}^{0, (q-2)}) \neq 0\}$ by exchanging the past and future parts. Let us choose the vectors $v_1^{0, (0)}, ..., v_{q-2}^{0, (q-3)}$. The choice of $v_{q-1}^{0, (q-2)}$ determines three connected components. If the parity of the number of future vectors among the timelike vectors of the $v_1^{0, (0)}, ..., v_{q-2}^{0, (q-3)}$ is unchanged by the involution, the component where $v_{q-1}^{0, (q-2)}$ is spacelike must be sent to itself and the two components where $v_{q-1}^{0, (q-2)}$ is past and future must be sent to each other. Since $i$ is an involution, the fourth connected component must be stable by $i$. If the parity of the number of future vectors among the timelike vectors of the $v_1^{0, (0)}, ..., v_{q-2}^{0, (q-3)}$ is changed by the involution, the component where $v_{q-1}^{0, (q-2)}$ is spacelike must not be stable by the involution while the two components defined by $v_{q-1}^{0, (q-2)}$ past and future are. This tells us that half of the non-positive components must be stable by $i$ while the other must be unstable, hence the result.

        Assume that $q$ is odd. The involution thus acts on $\{Q(v_{q-1}^{0, (q-2)}) \neq 0\}$ by preserving the past and future parts. If the parity of the number of future vectors among the timelike vectors of the $v_1^{0, (0)}, ..., v_{q-2}^{0, (q-3)}$ is unchanged by the involution, the three components given by the choice of $v_{q-1}^{0, (q-2)}$ must be preserved by $i$. If the parity of the number of future vectors among the timelike vectors of the $v_1^{0, (0)}, ..., v_{q-2}^{0, (q-3)}$ is changed by the involution, the components given by the choice of $v_{q-1}^{0, (q-2)}$ must not be stable by $i$. This gives us the result.
    \end{proof}

    \begin{rem}\label{componentwithspheres}
        In particular, there always exists at least one connected component which is stable by involution. The component where all the $v_j^{0, (j-1)}$ are spacelike is of particular interest as it is easy to exhibit antipodal spheres of higher dimensions contained inside. Let us take all $a_j^{i}, b_j^{i}$ equal to zero and embed $\R^{p-q+1,0}$ in $\R^{p-q+1,1}$ as a positive subspace. The map 
        \[\varphi : v \in \R^{p-q+1} \longmapsto \exp(u_1(v)) ... \exp(u_{q-1}(v)) \cdot F_0\]
        satisfies that for each distinct pair $v, v'$, the elements $\varphi(v)$ and $\varphi(v')$ are transverse. Adding $F_{\infty}$ as the point at infinity gives us an antipodal sphere of dimension $p-q+1$. It is not known to us whether this sphere may be realized as the boundary of a $P_{1,...,q-1}$-Anosov subgroup, except when $q=2$ where those are the boundaries of fuchsian representations of uniform lattices of $\SO_0(p,1)$.
    \end{rem}

    \begin{rem}
        When considering the space of full flags $\SO_0(p,q)/P_{1,...,q}$ for $p > q+1$, since the last equation does not change the count of the connected components, the computation of the involution remains the same.
    \end{rem}

    \subsection{\texorpdfstring{$q$}{q}-flags in \texorpdfstring{$\SO_0(q+1,q)$}{SO(q+1,q)}}

    Let us compute the involution for the space of maximal photons $\SO_0(q+1,q)/P_q$. The space of full flags will be discussed at the end of the section.

    Dealing only with the maximal photons is formally equivalent to taking all the $a_j^{i}$ equal to zero, which brings us back to Proposition \ref{involutioncompute}. The equation of transversality is given by $f^2 \neq 0$, where $f$ is defined in remark \ref{functionf} and the two connected components of $\Omega(F_0) \cap \Omega(F_{\infty})$ are determined by the sign of $f$. When all necessary vectors are non-zero, $f$ is equal to $v_1^{0, (0)}...v_q^{0, (q-1)}$. The same reasoning as in \ref{involutioncompute} tells us that depending on the parity of $j$, $(v_j^{0, (j-1)} < 0)$ and $(v_j^{0, (j-1)} > 0)$ are either preserved or sent to each other. One may easily determine in which case the sign of $f$ is reversed and in which case it is invariant.

    \begin{prop}
        When $q=1$ or $q=2$ mod $4$, the two connected components of $\SO_0(q+1,q)/P_q$ are swapped by $i$.
        When $q=3$ or $q=0$ mod $4$, the two connected components of $\SO_0(q+1,q)/P_q$ are stable by $i$.

    \end{prop}

    From this we may deduce the following theorem using Theorem \ref{DGR2}:

    \begin{teo}\label{th:rigid_q+1_q}
        Assume that $q = 1$ or $q=2$ mod $4$ and let $\Gamma$ be a $P_{\Theta}$-Anosov subgroup of $\SO_0(q+1,q)$ where $\Theta$ contains the last root. Then $\Gamma$ is virtually isomorphic to either a free group or a surface group.
    \end{teo}

    This applies in particular to Borel-Anosov subgroups of $\SO_0(q+1,q)$; however this result was already known due to Dey \cite{dey2024borelanosovsubgroupsrm}.

    \subsection{\texorpdfstring{$(q-1)$}{(q-1)}-flags in \texorpdfstring{$\SO_0(q,q)$}{SO(q,q)}}

    As stated in Proposition \ref{prop:countSOqq}, when $q$ is even, the space $\Omega(F_0) \cap \Omega(F_{\infty})$ is the union of two connected components. Since all $a_j^{i}$ and $v_j^{i}$ are equal to zero, the involution becomes $i(b_j^{i}) = -b_j^{i}$.

    \begin{prop}
        When $q = 0$ mod $4$, the two connected components of $\Omega(F_0) \cap \Omega(F_{\infty})$ are stable under the involution.
        When $q = 2$ mod $4$, the two connected components of $\Omega(F_0) \cap \Omega(F_{\infty})$ are swapped by the involution.

    \end{prop}

    \begin{proof}
        The two connected components are determined by the sign of the Pfaffian of $S$. The effect of the involution on $S$ is $i(S) = -S$. When the size of $S$ is equal to $0$ mod $4$, $\mathrm{Pf}(-S) = \mathrm{Pf}(S)$ while when it is equal to $2$ mod $4$, $\mathrm{Pf}(-S) = - \mathrm{Pf}(S)$, hence the result.
    \end{proof}

    Using Theorem \ref{DGR2}, we obtain the following result:

    \begin{teo}\label{th:rigid_q_q}
        Assume that $q = 2$ mod $4$ and let $\Gamma$ be a $P_{\Theta}$-Anosov subgroup of $\SO_0(q,q)$ with $\Theta$ containing one of the two last roots. Then $\Gamma$ is virtually isomorphic to either a free group or a surface group.
    \end{teo}

    \section{Constructing higher dimensional Anosov subgroups}

    This section will be dedicated to the construction of $P_{\Theta}$-Anosov subgroups of $\SO_0(p,q)$ which are not virtually isomorphic to a surface group or a free group for some of the $\Theta$ which do not fall under the assumptions of Theorems \ref{th:rigid_q_q} or \ref{th:rigid_q+1_q}. It will heavily rely on a combination theorem by Dey and Kapovich, see \cite{Dey_2023}, \cite{dey2024kleinmaskitcombinationtheoremanosov}. For $\Gamma$ a $P_{\Theta}$-Anosov subgroup, we will denote by $\partial \Gamma$ the image of the Gromov boundary of $\Gamma$ by the Anosov map associated to $\Gamma$ into the corresponding flag variety.

    \subsection{In \texorpdfstring{$(1,\dots,q-2)$}{(1,...,q-2)}-flags in \texorpdfstring{$\SO_0(q+1,q)$}{SO(q+1,q)}}

    Let us first see that, given a pair of transverse points, it is easy to build a $P_{\Theta}$-Anosov cyclic subgroup having those two points as boundary.
    
    \begin{lemma}
        Let $P_{\Theta} \subset \SO_0(p,q)$ be a parabolic subgroup and let $x^+, x^-$ be two transverse points in $\SO_0(p,q)/P_{\Theta}$. Then there exists $\delta \in \SO_0(p,q)$ such that $\langle \delta \rangle$ is $P_{\Theta}$-Anosov with boundary $\{x^+, x^-\}$.
    \end{lemma}

    \begin{proof}
        Since $x^+$ and $x^-$ are transverse, they can be completed into two transverse complete flags $\overline{x^+}, \overline{x^-} \in \SO_0(p,q)/P_{1, ..., q}$. We thus have $\overline{x^{\pm}} = (F^{\pm}_1, ..., F^{\pm}_q)$ with $F^-_i \oplus F^+_i \simeq \R^{i,i}$ for each $i$. Let us decompose orthogonally $F^-_q \oplus F^+_q = \R_1^{1,1} \oplus ... \oplus \R_q^{1,1}$ such that for each $i$, $F^-_i \oplus F^+_i = \R_1^{1,1} \oplus ... \oplus \R_i^{1,1}$ and let $V = (F^-_q \oplus F^+_q)^{\perp}$. Finally, let $\delta$ be the transformation acting trivially on $V$ and 

        \[\delta|_{\R_i^{1,1}} = \left(\begin{array}{cc}
         \lambda^i & 0\\
         0 & \lambda^{-i}\\
    \end{array}\right),\]

    with $\lambda > 1$. The subgroup $\langle \delta \rangle$ is $P_{1,...,q}$-Anosov with boundary $\{\overline{x^+}, \overline{x^-}\}$ in $\SO_0(p,q)/P_{1,...,q}$. In particular, $\langle \delta \rangle$ is also $P_{\Theta}$-Anosov with boundary $\{x^+, x^-\}$, hence the result.
    \end{proof}

    \begin{prop}\label{constructionanosov}
        Let $P_{\Theta}$ be a parabolic subgroup of $\SO_0(p,q)$, $\Gamma$ a $P_{\Theta}$-Anosov subgroup of $\SO_0(p,q)$ and $x^+, x^-$ two points of $\SO_0(p,q)/P_{\Theta}$ which are transverse to each other and such that for all $y \in \partial \Gamma \subset \SO_0(p,q)/P_{\Theta}$, $x^{\pm}$ and $y$ are transverse. Let $\delta \in \SO_0(p,q)$ associated to $x^+$ and $x^-$ by the previous lemma. Then there exists a finite index subgroup $\Gamma' \subset \Gamma$ and $k \in \mathbb{N}$ such that $\langle \Gamma', \delta^k \rangle$ is $P_{\Theta}$-Anosov and isomorphic to the free product $\Gamma' \ast\mathbb{Z}$.
    \end{prop}

    \begin{proof}
        This is a direct application of a result by Dey and Kapovich \cite{Dey_2023}, \cite{dey2024kleinmaskitcombinationtheoremanosov}. Since $x^+$ and $x^-$ are transverse to the whole boundary of $\Gamma$, there must exist pairwise transverse compacts $B$, $A^+$, $A^-$ such that $\partial \Gamma \subset B^{\circ}$ and $x^{\pm} \in (A^{\pm})^{\circ}$. For each $\gamma \in \Gamma$ let $|\gamma|$ be the word length associated to a finite set of generators. Since $\Gamma$ is $P_{\Theta}$-Anosov, there exists $K_1, K_2 >0$ such that for each $i \in \Theta$, 
        
        \[K_1 |\gamma| - K_2 < \log \left( \frac{\sigma_i(\gamma)}{\sigma_{i+1}(\gamma)} \right).\] 
        
        In particular, the sequence $(\gamma^n)$ has an attractive point $\gamma^+$ in $\SO_0(p,q)/P_{\Theta}$, and since $A^{\pm}$ are both transverse to $\gamma^+ \in \partial \Gamma$, there exists a constant $C > 0$ such that if $|\gamma| > C$ then $\gamma \cdot A^{\pm} \subset B^{\circ}$.
        
         From this we deduce that there exists a finite index subset $\Gamma' \subset \Gamma$ such that for all non trivial $\gamma \in \Gamma'$, $\gamma \cdot A^{\pm} \subset B$. Inversely, there exists $k$ such that $\delta^k \cdot B \subset A^+$ and $\delta^{-k} \cdot B \subset A^-$. By applying the result from Dey and Kapovich, we then get that $\langle \Gamma', \delta^k \rangle$ is $P_{\Theta}$-Anosov and isomorphic to $\Gamma' \ast \mathbb{Z}$.
    \end{proof}

    \begin{rem}
        When $\Gamma$ is a surface group, any finite index subgroup $\Gamma' \subset \Gamma$ is also a surface group.
    \end{rem}
    
    \begin{rem}
        It is actually not necessary to take a finite index subgroup of $\Gamma$. For each $\gamma \neq e$ in $\Gamma$, the set $O_{\gamma}$ of points $x$ in $\SO_0(q+1,q)/P_{\Theta}$ which are transverse to $\gamma \cdot x$ is an open and dense subset of $\SO_0(q+1,q)/P_{\Theta}$. Since $\Gamma$ is countable, the set $\bigcap_{\gamma \neq e} O_{\gamma}$ is still dense, meaning that up to a small perturbation of $x^+$ and $x^-$ one may assume that they are both in $\bigcap_{\gamma \neq e} O_{\gamma}$. We may then take $B$ to be the reunion of a small compact containing $\partial \Gamma$ and of small enough compacts containing the $\gamma \cdot x^{\pm}$. There is only a finite number of those compacts which are not contained within the compact containing $\partial \Gamma$, thus $B$ is also compact. We then get that $A^{\pm}$ and $B$ are transverse to each other, hence the result.
    \end{rem}

    Since the free product of a surface group with $\mathbb{Z}$ is not virtually isomorphic to either a surface group or a free group, this means that in order to construct a suitable $P_{\Theta}$-Anosov subgroup one only needs to finds a $P_{\Theta}$-Anosov subgroup $\Gamma$ isomorphic to a surface group and two transverse points $x^+$ and $x^-$ in $\SO_0(p,q)/P_{\Theta}$ which are transverse to the boundary of $\Gamma$.

    \begin{prop}
        Let $P_{\Theta} = P_{1, ..., q-2}$ be the stabilizer of a $(1,\dots,q-2)$-flag. There exists a $P_{\Theta}$-Anosov subgroup of $\SO_0(p,q)$ which is isomorphic to the free product of a surface group and an infinite cyclic group.
    \end{prop}

    \begin{proof}
        Let $\Gamma$ be a $\Theta$-positive subgroup of $\SO_0(q,q-1)$. The inclusions $$\Gamma \subset \SO_0(q,q-1) \subset \SO_0(q+1, q)$$ gives a $P_{1, ..., q-2}$-Anosov subgroup of $\SO_0(q+1,q)$. For each $k \leq q-2$, let $$\pi_k : \SO_0(q+1,q)/P_{1,...,q-2} \rightarrow \SO_0(q+1,q)/P_k$$ be the projection on the $k$-th factor. Let $x_0$ and $x_{\infty}$ be two distinct points in $\partial \Gamma$. In the affine chart of $\Ein_{q,q-1}$ defined by $\pi_1(\xi(x_0))$ and $\pi_1(\xi(x_{\infty}))$, the image of $\pi_1 \circ \xi$ is contained within a subspace $V$ of $\R^{q, q-1}$ of signature $(q-1, q-2)$. Furthermore, depending on whether $q$ is even or odd, all the elements of the image are of positive or negative norm. Let $x^+_1$ and $x^-_1$ be two elements of the orthogonal of $V$ in the affine chart such that $x_1^- = - x_1^+$ and with norm of the same sign as the image of $\pi_1 \circ \xi$. For each $x \in \partial \Gamma$, $x^{\pm}_1$ is transverse to $\pi_1(\xi(x))$ as transversality between two points in an affine chart is equivalent the segment between the two points not being lightlike.

        Let us consider the set of flags of affine isotropic spaces in $\R^{q, q-1}$ of dimension $1, ..., q-3$ containing $x_1^{\pm}$. Since the intersection $E$ of the isotropic cones of $x_1^{\pm}$ and $\xi(x_{\infty})$ is a copy of $\Ein_{q-1, q-2}$, this space can be identified with $\SO_0(q+1,q)/P_{1,...,q-3}$, and since $x_1^{\pm}$ is orthogonal to $V$, the boundary of $V$ in the isotropic cone of $\xi(x_{\infty})$ is contained within $E$. We can then do the same procedure in the affine chart of $E$ defined by $\pi_2(\xi(x_{\infty}))$ and $\pi_2(\xi(x_0))$ to get two flags $(x_1^{\pm},...,x_{q-2}^{\pm})$ which are transverse to each other and transverse to the image of $\xi$. This yields the result via Proposition \ref{constructionanosov}. Furthermore, the inclusion $\SO_0(q+1,q) \subset \SO_0(p,q)$ extends the result for $p \geq q+1$.
    \end{proof}

    \begin{rem}
        Iterating this procedure allows one to construct a $P_{1,...,q-2}$-Anosov subgroup which is isomorphic to $\Gamma \ast F_n$ for each $n$. Indeed, theorem A from \cite{Dey_2023} tells us that the boundary of $\Gamma \ast \mathbb{Z}$ is contained within $B \cup A^{\pm}$; by taking $A^{\pm}$ small enough, it is possible to find another pair of transverse points transverse to $B \cup A^{\pm}$ and thus to the boundary of $\Gamma \ast \mathbb{Z}$. The result follows by repeating the process for $\Gamma \ast F_{n-1}$.
    \end{rem}

    \subsection{In \texorpdfstring{$(1,\dots,q)$}{(1,...,q)}-flags in \texorpdfstring{$\SO_0(3q+1,q)$}{SO(3q+1,q)}}

    Using the same methods for $\Theta = \{1,...,q\}$ in $\SO_0(p,q)$, it is enough to find a $P_{\Theta}$-Anosov surface group $\Gamma$ and two points $x^{\pm}$ in $\SO_0(p,q)/P_{\Theta}$ which are transverse to $\partial \Gamma$.

    \begin{lemma}\label{decreasingreps}
        There exists quasi-fuchsian representations $\rho_1,...,\rho_q$ from $\Gamma$ to $\SO_0(3,1)$ such that for each $i$, $\rho_i$ strictly dominates $\rho_{i+1}$.
    \end{lemma}

    \begin{proof}
        We know from Deroin--Tholozan (\cite{TholozanDeroin}, theorem A) that for each quasi-fuchsian representation $\rho : \Gamma \to \SO_0(3,1)$ which is not fuchsian, there exists a fuchsian representation $\overline{\rho}$ which strictly dominates $\rho$. Let $\rho_q$ be a quasi-fuchsian representation which is not fuchsian and let $\overline{\rho_q}$ be a fuchsian representations strictly dominating $\rho_q$. Strictly dominating a fixed representation is an open condition; one way to see it is to use claim 2.2 in \cite{BG04} to show that for each $\varepsilon > 0$ there exists a neighborhood $V$ of $\overline{\rho_q}$ in the caracter variety sur that for each representation $\rho$ in $V$ there is a $(1+\varepsilon)$-bilipschitz diffeomorphism between the convex cores of $\rho$ and $\overline{\rho_q}$. Combined with the fact that $\overline{\rho_q}$ dominates $\rho_q$, we get that any $\rho$ in $V$ dominates $\rho_q$ as well. We may then slightly deform $\overline{\rho_q}$ into a quasi-fuchsian representation $\rho_{q-1}$ which strictly dominates $\rho_q$ and is not fuchsian. Repeating this process gives the result.
    \end{proof}

    While building $P_{\Theta}$-Anosov representations of $\Gamma$ into $\SO_0(p,q)$ is easy by taking embeddings of Hitchin representations of $\SO_0(q+1,q)$, it is not possible to find a pair of transverse points in $\SO_0(p,q)/P_{\Theta}$ which are transverse to $\partial \Gamma$. Using the previous lemma allows us to build another sort of $P_{\Theta}$-Anosov representation of $\Gamma$.

    \begin{lemma}
        Let $(\rho_1,...,\rho_q)$ be as in \ref{decreasingreps}. Let $\rho$ be the representation from $\Gamma$ to $\SO_0(3q,q)$ obtained by diagonal inclusion of the $\rho_i$ in $\SO_0(3,1)$ in $\SO_0(3q,q)$. Then $\rho$ is $P_{\Theta}$-Anosov.
    \end{lemma}

    \begin{proof}
        The fact that $\rho$ is $P_{\Theta}$-divergent follows from the fact that for each $i$, $\rho_i$ dominates $\rho_{i+1}$. The diagonal inclusion of $\SO_0(3,1)$ in $\SO_0(3q,q)$ is equivalent to an ordered orthogonal splitting $\R^{3q,q} = \R_1^{3,1} \oplus^{\perp} ... \oplus^{\perp} \R_q^{3,1}$ and taking an isotropic vector in each $\R_i^{3,1}$ gives a complete flag of isotropic spaces in $\R^{3q,q}$. Since the space of isotropic vectors in each $\R_i^{3,1}$ is a spacelike two-dimensional sphere, the boundary $\partial \Gamma$ must be totally transverse in the space of complete flags of photons. In fact, for two distinct points $x_{\infty}, x_0$ in $\partial \Gamma$, each segment $(x_0, x_{\infty})$ in $\partial \Gamma$ lies in the connected component described in \ref{componentwithspheres}.
    \end{proof}

    Let us now build $\rho : \Gamma \rightarrow \SO_0(3q+1,q)$ by combining the previous representation in $\SO_0(3q,q)$ with the inclusion $\SO_0(3q,q) \subset \SO_0(3q+1,q)$. We will proceed with the same method as in the previous section to find two transverse points which are transverse to $\partial \Gamma$. For each $k \leq q$, let $$\pi_k : \SO_0(3q+1,q)/P_{\Theta} \rightarrow \SO_0(3q+1,q)/P_k$$ be the projection on the $k$-th factor. Let $x_0$ and $x_{\infty}$ be two distinct points in $\partial \Gamma$. In the affine chart of $\Ein_{3q+1,q}$ defined by $\pi_1(\xi(x_0))$ and $\pi_1(\xi(x_{\infty}))$, the image of $\pi_1 \circ \xi$ is contained within a subspace $V$ of $\R^{3q+1, q}$ of signature $(3q, q)$. Furthermore, all the elements of the image are of positive norm. Let $x^+_1$ and $x^-_1$ be two elements of the orthogonal of $V$ such that $x_1^- = - x_1^+$; such elements must be of positive norm. For each $x \in \partial \Gamma$, $x^{\pm}_1$ is transverse to $\pi_1(\xi(x))$ as the segment joining the two must be positive.

        Let us consider the set of flags of affine isotropic spaces in $\R^{3q+1, q}$ of dimension $1, ..., q-1$ containing $x^{\pm}$. Since the intersection $E$ of the isotropic cones of $x^{\pm}$ and $x_{\infty}$ is a copy of $\Ein_{3q, q-1}$, this space can be identified with $\SO_0(3q,q-1)/P_{1,...,q-1}$, and since $x^{\pm}$ is orthogonal to $V$, the boundary of $V$ in the isotropic cone of $x_{\infty}$ is contained within $E$. We can then do the same procedure in the affine chart of $E$ defined by $\pi_2(\xi(x_{\infty}))$ and $\pi_2(\xi(x_0))$ to get two flags $(x_1^{\pm},...,x_{q}^{\pm})$ which are transverse to each other and transverse to the image of $\xi$. This yields the result via Proposition \ref{constructionanosov}. Furthermore, the inclusion $\SO_0(3q+1,q) \subset \SO_0(p,q)$ extends the result for $p \geq 3q+1$.

    \printbibliography

 \end{document}